\documentclass{amsart}  
\usepackage{amscd, amssymb}
\newtheorem{thm}{Theorem}[section] % Theorem style, bold title, italic text
\newtheorem{cor}[thm]{Corollary}
\newtheorem{lem}[thm]{Lemma}
\newtheorem{prop}[thm]{Proposition}

\theoremstyle{definition} % Definition style - Bold title, roman text
\newtheorem{dfn}[thm]{Definition}
\newtheorem{dfns}[thm]{Definitions}
\newtheorem{example}[thm]{Example}
\newtheorem{examples}[thm]{Examples}

\theoremstyle{remark} % Remark style - Italised title, roman text
\newtheorem{rmk}[thm]{Remark}
\newtheorem{rmks}[thm]{Remarks}

\begin{document}

\title{Higher Rank Graph $C^*$-algebras}

\author[Kumjian]{Alex Kumjian}
\address{Department of Mathematics (084)\\
University of Nevada\\
Reno NV 89557--0045, USA.}
\email{alex@unr.edu}
\thanks{Research of the first author partially supported by {\sc nsf} grant
DMS-9706982} 
\author[Pask]{David Pask}
\address{Department of Mathematics\\
University of Newcastle\\
NSW 2308, Australia}
\email{davidp@maths.newcastle.edu.au}
\thanks{Research supported by University of Newcastle RMC project grant}

\date{Jan.\ 2000: {\em New York J.\ Math.} {\bf 6} (2000), 1-20.
           {\tt http://nyjm.albany.edu:8000/j/2000/6-1nf.htm}}

\newcommand{\abstracttext}{}

\begin{abstract}
\abstracttext
Building on recent work of  Robertson and Steger, we
associate a $C^*$--algebra to a combinatorial object which may be thought
of as a higher rank graph.  This $C^*$--algebra is
shown to be isomorphic to that of the associated path groupoid.   Sufficient
conditions on the higher rank graph are found for the associated
$C^*$--algebra to be simple, purely infinite and AF.   Results concerning the
structure of crossed products by certain natural actions of discrete groups
are obtained; a technique for constructing rank $2$ graphs from ``commuting''
rank $1$ graphs is given.
\end{abstract}

\maketitle
\tableofcontents

%%%% **** The text of the paper starts here **** %%%%

%\section{Introduction}
\noindent
In this paper we shall introduce the notion of a higher rank graph and
associate a $C^*$--algebra to it in such a way as to generalise the
construction of the $C^*$--algebra of a directed graph as studied in
\cite{ck,kprr,kpr} (amongst others). Graph $C^*$--algebras include up to
strong Morita equivalence Cuntz--Krieger algebras and AF algebras. The
motivation for the form of our generalisation comes from the recent work of
Robertson and Steger \cite{rs1,rs2,rs3}. In \cite{rs1} the authors study
crossed product $C^*$--algebras arising from certain
group actions on $\tilde{A_2}$-buildings and show that they are generated
by two families of partial isometries which satisfy certain relations amongst
which are Cuntz--Krieger type relations
\cite[Equations (2), (5)]{rs1} as well as more intriguing
commutation relations \cite[Equation (7)]{rs1}. In \cite{rs2} they give
a more general framework for studying such algebras involving certain
families of commuting $0-1$ matrices. In particular the associated
$C^*$--algebras are simple, purely infinite and generated by a family of
Cuntz--Krieger algebras associated to these matrices. It is this
framework which we seek to cast in graphical terms to include a wider
class of examples (including graph $C^*$--algebras). 

What follows is a brief outline of the paper. In the first section we
introduce the notion of a higher rank graph as a purely combinatorial
object: a small category $\Lambda$ gifted with a degree map $d :\Lambda
\rightarrow {\bf N}^k$ (called shape in \cite{rs2}) playing the role of
the length function. No detailed knowledge of category theory is
required to read this paper. The associated $C^*$--algebra $C^* ( \Lambda )$
is defined as the universal $C^*$--algebra
generated by a family of partial isometries 
$\{ s_\lambda : \lambda \in \Lambda \}$ 
satisfying relations similar to those of
\cite{kpr} (our standing assumption is that our higher rank graphs
satisfy conditions analogous to a directed graph being row--finite and
having no sinks). We then describe some basic examples and indicate the
relationship between our formalism and that of \cite{rs2}.

In the second section we introduce the path groupoid ${\mathcal G}_\Lambda$
associated to a higher rank graph $\Lambda$ (cf.~\cite{r,d,kprr}).
Once the infinite path space $\Lambda^\infty$ is formed (and a few
elementary facts are obtained) the construction is fairly routine.
It follows from  the gauge-invariant uniqueness theorem (Theorem
\ref{giuniqueness}) that 
$C^* ( \Lambda) \cong C^* ({\mathcal G}_\Lambda )$.  
By the universal property $C^* (\Lambda )$ carries a canonical action
of ${\bf T}^k$ defined by
\begin{equation} \label{gaugedef}
\alpha_t ( s_\lambda ) = t^{d ( \lambda )} s_\lambda 
\end{equation}
\noindent
called the gauge action. In the third section we prove the gauge--invariant
uniqueness theorem, which is the key result for analysing $C^*( \Lambda )$
(cf.~\cite{bprs,ahr}, see also \cite{ck,rs2} where similar techniques are used
to prove simplicity). It gives conditions under which a homomorphism with domain
$C^* ( \Lambda )$ is faithful: roughly speaking, if the homomorphism
is equivariant 
for the gauge action and nonzero on the generators then it is faithful.
This theorem has a number of interesting consequences, amongst which are the
isomorphism mentioned above and the fact that the higher rank
Cuntz--Krieger algebras of \cite{rs2} are isomorphic to
$C^*$--algebras associated to suitably chosen higher rank graphs.

In the fourth section we characterise, in terms of an aperiodicity
condition on $\Lambda$, the circumstances under which the groupoid
${\mathcal G}_\Lambda$ is essentially free. This aperiodicity condition
allows us to prove a second uniqueness theorem analogous to the
original theorem of \cite{ck}.  In \ref{simple} and \ref{pi} we 
obtain conditions under which $C^*( \Lambda )$ is simple and purely infinite 
respectively which are similar to those in \cite{kpr} but with the 
aperiodicity condition replacing condition (L). 
%(in the sense that every hereditary subalgebra contains an infinite
%projection) 

In the next section we show that, given a functor $c : \Lambda
\rightarrow G$ where $G$ is a discrete group, then as in \cite{kp} one may
construct a skew product $G \times_c \Lambda$ which is also a higher
rank graph. If $G$ is abelian then there is a natural action
$\alpha^c : \widehat{G} \rightarrow \mbox{Aut} \, C^* ( \Lambda )$
such that
\begin{equation} \label{dualdef}
\alpha^c_\chi ( s_\lambda ) = \langle \chi , c (  \lambda ) \rangle
s_\lambda ;
\end{equation}
\noindent
moreover 
$C^* ( \Lambda ) \rtimes_{\alpha^c} \widehat{G} \cong 
C^*(G  \times_c  \Lambda )$. 
Comparing (\ref{gaugedef}) and (\ref{dualdef}) we see that the gauge action
$\alpha$ is of the form $\alpha^d$ and as a consequence we may show that the
crossed product of $C^* ( \Lambda )$ by the gauge action is
isomorphic to $C^* ( {\bf Z}^k \times_d \Lambda )$; this $C^*$--algebra is then
shown to be AF.  By Takai duality $C^* ( \Lambda )$ is strongly Morita
equivalent to a crossed product of this AF algebra by the dual action of
${\bf Z}^k$.  Hence $C^* ( \Lambda )$ belongs to the bootstrap class 
${\mathcal N}$ of $C^*$--algebras for which the UCT applies 
(see \cite{rs-uct}) and is consequently nuclear.  If a discrete group
$G$ acts freely on a $k$-graph $\Lambda$, then the quotient object
$\Lambda / G$ inherits the structure of a $k$--graph; moreover 
(as a generalisation of \cite[Theorem 2.2.2]{gt}) there is a functor 
$c :\Lambda / G \rightarrow G$ such that 
$\Lambda \cong G \times_c ( \Lambda/ G )$ in an equivariant way. This
fact allows us to prove that
\[
C^* ( \Lambda ) \rtimes G \cong C^* ( \Lambda / G ) \otimes 
{\mathcal K} \left(\ell^2 ( G ) \right)
\]
where the action of $G$ on $C^*(\Lambda )$ is induced from that on $\Lambda$.
Finally in \S 6 a technique for constructing a $2$-graph from ``commuting''
$1$-graphs  $A, B$ with the same vertex set is given.  The
construction depends on 
the choice of a certain bijection  between pairs of composable edges:
$\theta :(a,b) \mapsto (b', a')$ where $a, a' \in A^1$ and 
$b, b' \in B^1$; the resulting $2$-graph is denoted $A *_\theta B$.
It is not hard to show that every $2$-graph is of this form.

Throughout this paper we let ${\bf N} = \{ 0 , 1 , \ldots \}$ denote the
monoid of natural numbers under addition. 
For $k \ge 1$  regard ${\bf N}^k$ as an abelian monoid under addition
with identity $0$ (it will sometimes be useful to regard ${\bf N}^k$
as a small category with one object) and canonical generators $e_i$ for 
$i=1 , \ldots, k$; we shall also regard ${\bf N}^k$ as 
the positive cone of ${\bf Z}^k$ under the usual coordinatewise partial order:
thus $m \leq n$ if and only if $m_i \le n_i$ for all $i$ where 
$m = ( m_1 , \ldots , m_k )$, and $n = ( n_1 , \ldots , n_k )$ (this
makes ${\bf N}^k$ a lattice).  

We wish to thank Guyan Robertson and Tim Steger for providing
us with an early version of their paper \cite{rs2}; the first author would also
like to thank them for a number of stimulating conversations and the staff of 
the Mathematics Department at Newcastle University for their hospitality 
during a recent visit.

\section{Higher rank graph $C^*$--algebras}

\noindent
In this section we first introduce what we shall call a higher rank graph
as a purely combinatorial object (we do not know whether this concept
has been studied before). Our definition of a higher rank
graph is modelled on the path category of a directed graph (see
\cite{h}, \cite{mu}, \cite[\S II.7]{macl} and Example \ref{firstex}). 
Thus a higher rank graph will be defined to be a small category
gifted with a degree map (called shape in \cite{rs2}) satisfying a
certain factorization property.
We then introduce the associated $C^*$--algebra whose definition is
modelled on that of the $C^*$--algebra of a graph as well as the
definition of \cite{rs2}.

\begin{dfns} \label{rkgdef}
A {\bf $k$-graph} (rank $k$ graph or higher rank graph) $(\Lambda, d)$
consists of a countable small category $\Lambda$ (with range and 
source maps $r$ and $s$ respectively) together with a functor 
$d : \Lambda \rightarrow {\bf N}^k$ satisfying the {\bf  factorisation 
property:} for every $\lambda \in \Lambda$ and $m, n \in {\bf
N}^k$ with $d( \lambda ) = m + n$, there are unique elements
$\mu , \nu \in \Lambda$ such that $\lambda = \mu \nu $ and
$d( \mu ) = m$, $d( \nu ) = n$. For $n \in {\bf N}^k$ we write
$\Lambda^n := d^{-1} (n)$. A morphism between $k$-graphs 
$(\Lambda_1 , d_1 )$ and $( \Lambda_2 , d_2 )$ is a functor 
$f : \Lambda_1 \rightarrow \Lambda_2$ compatible with the degree maps.
% that is $d_2 ( f ( \lambda ) ) =  d_1 ( \lambda )$ for all 
% $\lambda \in \Lambda_1$. 
\end{dfns}

\begin{rmks}
The factorisation property of \ref{rkgdef} allows us to identify
$\mbox{Obj} ( \Lambda )$, the objects of $\Lambda$ with 
$\Lambda^0$. Suppose $\lambda \alpha = \mu \alpha$ in $\Lambda$
then by the the factorisation property $\lambda = \mu$; left
cancellation follows similarly. We shall write the objects of $\Lambda$
as $u, v, w , \ldots$ and the morphisms as greek letters $\lambda ,
\mu , \nu \ldots$.  We shall frequently refer to $\Lambda$ as a
$k$-graph without mentioning $d$ explicitly. 
\end{rmks}

\noindent
It might be interesting to replace ${\bf N}^k$ in Definition \ref{rkgdef}
above by a monoid or perhaps the positive cone of an ordered abelian group.
% but we have no meaningful examples at hand to motivate this extension.

Recall that $\lambda , \mu \in  \Lambda$ are composable if and
only if $r ( \mu ) = s ( \lambda )$, and then $\lambda \mu \in
\Lambda$; on the other hand two finite paths $\lambda , \mu$ in
a directed graph may be composed to give the path $\lambda
\mu$ provided that $r ( \lambda ) = s ( \mu )$; so in \ref{firstex} 
below we will need to switch the range and source maps. 

\begin{example} \label{firstex}
Given a $1$-graph $\Lambda$,  define $E^0 = \Lambda^0$ and 
$E^1 = \Lambda^1$. 
If we define $s_E ( \lambda) = r ( \lambda )$ and 
$r_E(\lambda ) = s ( \lambda )$ then the quadruple 
$( E^0, E^1, r_E, s_E )$ is a directed graph in the sense of
\cite{kpr,kp}. On the other hand, 
given a directed graph $E = ( E^0 , E^1 , r_E , s_E )$, then 
$E^* = \cup_{n \ge 0} E^n$, the collection of finite paths, may be
viewed as small category with range and source maps given by 
$s(\lambda ) =  r_E ( \lambda )$ and $r ( \lambda ) = s_E ( \lambda )$.
If we let $d : E^*  \rightarrow {\bf N}$ be the length function (i.e.\ 
$d ( \lambda ) = n$ iff $\lambda \in E^n$) then $( E^* , d)$ is
a $1$-graph. 
\end{example}

\noindent
We shall associate a $C^*$--algebra to a $k$-graph  in such a way that
for $k=1$ the associated $C^*$--algebra is the same as that of the
directed graph. We shall consider other examples later. 

\begin{dfns}
The $k$-graph $\Lambda$ is {\bf row finite} if
for each $m \in {\bf N}^k$ and $v \in \Lambda^0$ the set
$\Lambda^m (v) := \{ \lambda \in \Lambda^m : r ( \lambda ) = v
\}$ is finite. Similarly $\Lambda$ has {\bf no sources} if
$\Lambda^m (v) \neq \emptyset$ for all $v \in \Lambda^0$ and
$m \in {\bf N}^k$.
\end{dfns}

\noindent
Clearly if $E$ is a directed graph then $E$ is row finite (resp.\ has no
sinks) if and only if $E^*$ is row finite (resp.\ has no sources).
Throughout this paper we will assume (unless otherwise stated)
that any $k$-graph $\Lambda$ is row finite and has no sources,
that is 
\begin{equation} \label{stand}
0 < \# \Lambda^n (v) < \infty  ~\mbox{for every}~ v \in \Lambda^0 ~
\mbox{and} ~ n \in {\bf N}^k .
\end{equation}

\noindent
The Cuntz--Krieger relations \cite[p.253]{ck} and the relations given in 
\cite[\S 1]{kpr} may be interpreted as providing a representation of a
certain directed graph by partial isometries and orthogonal
projections.  This view motivates the definition of $C^* ( \Lambda )$.

\begin{dfns} \label{cstarlambdadef}
Let $\Lambda$ be a row finite $k$-graph with no sources. Then $C^* ( \Lambda )$ 
is defined to be the universal $C^*$--algebra generated by a family 
$\{ s_\lambda : \lambda \in \Lambda \}$ of partial isometries
satisfying:
\begin{itemize}
\item[(i)] $\{ s_v : v \in \Lambda^0 \}$ is a family of mutually 
orthogonal projections,
\item[(ii)] $s_{\lambda \mu} = s_\lambda s_\mu$ for all
$\lambda , \mu \in \Lambda$ such that $s ( \lambda ) = r ( \mu )$,
\item[(iii)] $s_\lambda^* s_\lambda = s_{s( \lambda )}$ for all
$\lambda \in \Lambda$,
\item[(iv)] for all $v \in \Lambda^0$ and $n \in {\bf N}^k$ we  have
$\displaystyle s_v = \sum_{\lambda \in
\Lambda^n ( v )} s_\lambda s_\lambda^*$.
\end{itemize}
\noindent
For $\lambda \in \Lambda$, define $p_\lambda = s_\lambda s_\lambda^*$
(note that $p_v = s_v$ for all $v \in \Lambda^0$).
A family of partial isometries satisfying (i)--(iv) above is called a 
{\bf $*$--representation} of $\Lambda$.
\end{dfns}

\begin{rmks} \label{rad}
\begin{itemize}
\item[(i)] If $\{ t_\lambda : \lambda \in \Lambda \}$ is a
$*$--representation of $\Lambda$ then the map $s_\lambda \mapsto
t_\lambda$ defines a $*$--homomorphism from  $C^* ( \Lambda )$ to 
$C^* ( \{ t_\lambda : \lambda \in \Lambda \} )$.
\item[(ii)]
If $E^*$ is the $1$-graph associated to the directed graph $E$ (see
\ref{firstex}), then by restricting a $*$--representation 
to $E^0$ and $E^1$ one obtains a Cuntz--Krieger family for $E$ in the
sense of \cite[\S 1]{kpr}.  
Conversely every Cuntz--Krieger family for $E$ extends uniquely to a 
$*$--representation of $E^*$.
\item[(iii)]
In fact we only need the relation (iv) above to be satisfied for
$n = e_i \in {\bf N}^k$ for $i=1 , \ldots , k$, 
the relations for all $n$ will then follow 
(cf.~\cite[Lemma 3.2]{rs2}). Note that the definition
of $C^* ( \Lambda )$ given in \ref{cstarlambdadef} may be extended
to the case where there are sources by only requiring that relation (iv)
hold for $n = e_i$ and then only if $\Lambda^{e_i} (v) \neq
\emptyset$ (cf.~ \cite[Equation (2)]{kpr}).
\item[(iv)]  For $\lambda , \mu \in \Lambda$ if 
$s ( \lambda ) \neq s (\mu )$ then $s_\lambda s_\mu^* = 0$.  The
converse follows from \ref{univexists}.
\item[(v)] Increasing finite sums of $p_v$'s form an approximate
identity for $C^* ( \Lambda )$ (if $\Lambda^0$ is finite then
$\sum_{v \in \Lambda^0} p_v$ is the unit for $C^* ( \Lambda )$). It
follows from relations (i) and (iv) above that for any $n \in {\bf N}^k$, 
$\{ p_\lambda : d ( \lambda ) = n \}$ forms a collection of orthogonal
projections (cf.~\cite[3.3]{rs2}); likewise increasing finite sums of these form an 
approximate identity for $C^* ( \Lambda )$ (see \ref{partition}).
\item[(vi)]  The above definition is not stated most efficiently.
Any family of operators
$\{ s_\lambda : \lambda \in \Lambda \}$ satisfying the above conditions must
consist of partial isometries.  The first two axioms could also be replaced by:
$$
s_\lambda s_\mu = 
\begin{cases}
s_{\lambda \mu} &\mbox{ if } s(\lambda) = r(\mu) \\
0 &  \mbox{ otherwise.} 
\end{cases}
$$
\end{itemize}
\end{rmks}

\begin{examples} \label{secondex}
\begin{itemize}
\item[(i)]  If $E$ is a directed graph, then by \ref{rad}(i) and (ii) we have
$C^* ( E^* ) \cong C^* (E)$ (see \ref{firstex}).  
\item[(ii)] For $k \ge 1$ let $\Omega = \Omega_k$ be the small
category with objects $\mbox{Obj} \, ( \Omega ) = {\bf N}^k$, and
morphisms 
$\Omega  = \{ (m,n) \in {\bf N}^k \times {\bf N}^k : m \le n \}$; 
the range and source maps are given by $r ( m , n ) = m$, $s ( m , n ) = n$. 
Let $d : \Omega\rightarrow {\bf N}^k$ be defined by $d ( m , n ) = n - m$. It is
then straightforward to show that $\Omega_k$ is a $k$-graph and
$C^* ( \Omega_k ) \cong {\mathcal K} \left( \ell^2 ( {\bf N}^k ) \right)$.
\item[(iii)] Let $T = T_k$ be the semigroup ${\bf N}^k$ viewed as a small 
category, then if $d : T \rightarrow {\bf N}^k$ is the identity map 
$( T , d )$ is a $k$-graph. It is not hard to show that 
$C^* ( T )\cong C ( {\bf T}^k )$, where $s_{e_i}$ for $1 \le i \le k$
are the canonical unitary generators.
\item[(iv)] 
Let $\{ M_1 , \ldots , M_k \}$ be square $\{0 ,1 \}$ matrices satisfying
conditions (H0)--(H3) of \cite{rs2} and let ${\mathcal A}$ be the
associated $C^*$-algebra. For $m \in {\bf N}^k$ let $W_m$ be the
collection of undecorated words in the finite alphabet $A$ of
shape $m$ as defined in \cite{rs2} then let
\[
W ~=~ \bigcup_{m \in {\bf N}^k} W_m .
\]
\noindent
Together with range and source maps $r ( \lambda ) = o ( \lambda )$, $s
( \lambda ) = t ( \lambda )$ and product defined in \cite[Definition
0.1]{rs2} $W$ is a small category. If we define $d
: W \rightarrow {\bf N}^k$ by $d ( \lambda ) =
\sigma ( \lambda )$, then one checks that $d$ satisfies the
factorisation property, and then from the second part of (H2) we
 see that $( W , d )$ is an
irreducible $k$-graph in the sense that for all $u,v \in
W_0$ there is $\lambda \in W$ such that $s ( \lambda )
= u$ and $r ( \lambda ) = v$.

We claim that the map $s_\lambda \mapsto s_{\lambda , s ( \lambda )}$
for $\lambda \in W$ extends to a *-homomorphism 
$C^* ( W ) \rightarrow {\mathcal A}$ for which 
$s_\lambda s_\mu^* \mapsto s_{\lambda ,\mu}$ (since these generate 
${\mathcal A}$ this will show that the map is onto).   
It suffices to verify that 
$\{ s_{\lambda , s ( \lambda )} : \lambda \in W \}$
constitutes a $*$--representation of $W$.  Conditions (i) and (iii)
are easy to check, (iv) follows from \cite[0.1c,3.2]{rs2}  with 
$u=v \in W^0$.  We check condition (ii):
if $s(\lambda) = r(\mu)$ apply \cite[3.2]{rs2} 
$$
s_{\lambda , s ( \lambda )}s_{\mu, s ( \mu )} =
\sum_{W^{d(\mu)}(s(\lambda))}s_{\lambda \nu, \nu }s_{\mu, s ( \mu )} =
s_{\lambda \mu, \mu }s_{\mu, s ( \mu )}  = s_{\lambda \mu, s ( \lambda\mu )} 
$$
where the sum simplifies using \cite[3.1, 3.3]{rs2} .
We shall show below that $C^* ( W ) \cong {\mathcal A}$.
\end{itemize}
\end{examples}

\noindent
We may combine higher rank graphs using the
following fact, whose proof is straightforward.

\begin{prop}
Let $( \Lambda_1 , d_1 )$ and $( \Lambda_2 , d_2 )$ be rank
$k_1$, $k_2$ graphs respectively, then $( \Lambda_1 \times
\Lambda_2 , d_1 \times d_2 )$ is a rank $k_1 + k_2 $ graph
where $\Lambda_1 \times \Lambda_2$ is the product category
and $d_1 \times d_2 : \Lambda_1 \times \Lambda_2
\rightarrow {\bf N}^{k_1 + k_2}$ is given by $d_1 \times d_2 (
\lambda_1 , \lambda_2 ) = ( d_1 ( \lambda_1 ) , d_2 ( \lambda_2 ) ) \in
{\bf N}^{k_1} \times {\bf N}^{k_2}$ for $\lambda_1 \in \Lambda_1$ and
$\lambda_2 \in \Lambda_2$.
\end{prop}

\noindent
An example of this construction is discussed in \cite[Remark
3.11]{rs2}. It is clear that   
$\Omega_{k+\ell} \cong \Omega_k \times \Omega_\ell$ 
where $k, \ell > 0$. 

\begin{dfn} \label{pullback}
Let $f : {\bf N}^\ell \rightarrow {\bf N}^k$ be a monoid
morphism, then if $( \Lambda ,d)$ is a  $k$-graph we may form the 
$\ell$-graph $f^* ( \Lambda )$ as follows:  (the objects of $f^*(
\Lambda )$ may be identified with those of $\Lambda$ and) $f^* (
\Lambda ) = \{ ( \lambda , n ) : d ( \lambda ) = f(n) \}$ with
$d ( \lambda , n ) = n$, $s ( \lambda , n ) = s ( \lambda )$ and 
$r ( \lambda , n ) = r ( \lambda )$. 
\end{dfn}

\begin{examples}\label{thirdex}
\begin{itemize}
\item[(i)] Let $\Lambda$ be a $k$--graph and put $\ell =1$, then if we
define the morphism $f_i (n) = n e_i$ for $1 \leq i \leq k$, we
obtain the {\bf coordinate graphs} $\Lambda_i : = f_i^* ( \Lambda )$
of $\Lambda$ (these are $1$--graphs).
\item[(ii)] Suppose $E$ is a directed graph and define 
$f : {\bf N}^2 \rightarrow {\bf N}$ by $( m_1, m_2) \mapsto m_1 + m_2$;
then the two coordinate graphs of $f^* ( E^* )$ are isomorphic to
$E^*$. We will show below that $C^* ( f^* ( E^* ) ) \cong 
C^* ( E^* ) \otimes C ( {\bf T} )$.
\item[(iii)] Suppose $E$ and $F$ are directed graphs and define
$f : {\bf N} \rightarrow {\bf N}^2$ by $f(m) = (m,m)$ then
$f^*(E^* \times F^*) = (E \times F)^*$ where  $E \times F$ denotes the
cartesian product graph (see \cite[Def.\ 2.1]{kp}).
\end{itemize}
\end{examples}

\begin{prop} \label{pifprop}
Let $\Lambda$ be a $k$-graph and $f : {\bf N}^\ell \rightarrow
{\bf N}^k $ a monoid morphism, then there is a $*$--homomorphism
$\pi_f : C^* ( f^* ( \Lambda ) ) \rightarrow C^* ( \Lambda )$ such
that  $s_{( \lambda , n )} \mapsto s_\lambda$; moreover if $f$ is
surjective, then $\pi_f$ is too.
\end{prop}

\begin{proof}
By \ref{rad}(i) it suffices to show that this is a $*$--representation of
$f^* ( \Lambda )$. Properties (i)--(iii) are straightforward to verify and
property (iv) follows by observing that for fixed $n \in {\bf N}^\ell$ and
$v \in \Lambda^0$ the map $f^* ( \Lambda )^n (v) 
\rightarrow \Lambda^{f(n)} (v)$ given by $( \lambda , n ) \mapsto
\lambda $ is a bijection. If $f$ is surjective, then it is clear that every
generator $s_\lambda$ of $C^* (
\Lambda )$ is in the range of $\pi_f$.
\end{proof}

\noindent
Later in \ref{giutconsequences} we will also show that $\pi_f$ is
injective if $f$ is injective.

\section{The path groupoid}

In this section we construct the path groupoid 
${\mathcal G}_{\Lambda}$ associated to a higher rank
graph $( \Lambda , d)$ along the lines of \cite[\S 2]{kprr}. 
%a second countable, locally compact, $r$--discrete (\'{e}tale)
Because some of the details are not quite the same as 
those in \cite[\S 2]{kprr} we feel it is useful to sketch the construction.
First we introduce the following analog of an infinite path in a higher rank graph:

\begin{dfns} \label{infpathdef}
Let $\Lambda$ be a $k$-graph, then 
\[
\Lambda^\infty ~=~ \{ x : \Omega_k  \rightarrow \Lambda  : x \;
\mbox{is a $k$-graph morphism} \} ,
\]
\noindent
is the infinite path space of $\Lambda$. For $v \in \Lambda^0$
let $\Lambda^\infty (v) = \{ x \in \Lambda^\infty : x ( 0 ) = v
\}$. For each $p \in {\bf N}^k$ define $\sigma^p : \Lambda^\infty
\rightarrow \Lambda^\infty$ by $\sigma^p (x) (m,n) = x(m+p,n+p)$
for $x \in \Lambda^\infty$ and $(m,n) \in \Omega$. (Note that
$\sigma^{p+q} = \sigma^p \circ \sigma^q$).
\end{dfns}

\noindent
By our standing assumption (\ref{stand}) one can show that for every $v
\in \Lambda^0$ we have $\Lambda^\infty (v)  \neq \emptyset$. Our
definition of $\Lambda^\infty$ is related to the definition of
$W_\infty$, the space of infinite words, given in the proof of \cite[Lemma
3.8]{rs2}.  If $E^*$ is the $1$-graph associated to the directed graph $E$ then
$(E^*)^\infty$ may be identified with $E^\infty$.

\begin{rmks} \label{pathdefine}
By the factorisation property the values of $x ( 0 , m )$ for 
$m \in {\bf N}^k$ completely determine $x \in \Lambda^\infty$.
To see this, suppose that $x ( 0 , m )$ is given for all $m \in
{\bf N}^k$ then for $( m , n ) \in  \Omega$,
$x ( m  , n )$ is the unique element $\lambda \in \Lambda$ 
such that $x ( 0 , n ) = x ( 0 , m ) \lambda$. 

More generally, let $\{ n_j : j \ge 0 \}$ be an increasing cofinal
sequence in ${\bf N}^k$ with $n_0 = 0$, then $x \in \Lambda^\infty$  
is completely determined by the values of $x ( 0 , n_j )$ (for example one
could take $n_j = j p$ where $p =  ( 1 , \ldots , 1 ) \in {\bf N}^k$).
Moreover, given a sequence $\{ \lambda_j : j \ge 1 \}$ in $\Lambda$ 
such that $s ( \lambda_j ) = r ( \lambda_{j+1} )$ and 
$d ( \lambda_j ) =  n_{j} - n_{j-1}$ there is a unique $x \in \Lambda^\infty$
such that $x ( n_{j-1} , n_j ) = \lambda_j$. For $(m,n) \in \Omega$ we
define $x(m,n)$ by the factorisation  property as follows: let $j$ be
the smallest index such that $n \le n_j$, then $x(m,n)$ is the unique
element of degree $n-m$ such that 
$\lambda_1\cdots\lambda_j =\mu x ( m,n ) \nu$ where $d ( \mu ) = m$
and $d ( \nu ) = n_j - n$.  It is straightforward to show that $x$ has
the desired properties.  
\end{rmks}

\noindent
We now establish a factorisation property for $\Lambda^\infty$ which is an easy
consequence of the above remarks: 

\begin{prop} \label{infpathfactor}
Let $\Lambda$ be a rank $k$ graph.  For all $\lambda \in \Lambda$ 
and $x \in \Lambda^\infty$ with $x(0) = s(\lambda)$, there is a unique 
$y \in \Lambda^\infty$ such that $x = \sigma^{d ( \lambda )} y$ and 
$\lambda = y ( 0,d (\lambda))$; we write $y = \lambda x$. 
Note that for every $x\in \Lambda^\infty$ and $p \in {\bf N}^k$ 
we have $x = x ( 0 , p) \sigma^p x$.
\end{prop}

\begin{proof}
Fix $\lambda \in \Lambda$ and $x \in \Lambda^\infty$ with 
$x(0) = s(\lambda)$. The sequence $\{n_j : j \ge 0 \}$ defined by 
$n_0 = 0$ and $n_j = (j-1)p + d(\lambda)$ for $j \ge 1$ is cofinal.  
Set $\lambda_1 = \lambda$ and $\lambda_j = x((j-2)p,(j-1)p)$ for $j \ge 2$
and let $y \in \Lambda^\infty$ be defined by the method given in 
\ref{pathdefine}.  Then $y$ has the desired properties.
\end{proof}

\noindent
Next we construct a  basis of compact open sets for the
topology on $\Lambda^\infty$ indexed by $\Lambda$:

\begin{dfns}
Let $\Lambda$ be a rank $k$ graph.
For $\lambda \in \Lambda$ define
\[
Z ( \lambda ) = \{ \lambda x \in \Lambda^\infty : s ( \lambda ) =
x (0)  \} = \{ x : x ( 0 , d ( \lambda ) ) = \lambda \} .
\]
\end{dfns}

\begin{rmks} \label{partition}
Note that $Z(v) = \Lambda^\infty(v)$ for all $v \in \Lambda^0$.
For fixed $n \in {\bf N}^k$ the sets 
$\{ Z ( \lambda ) : d ( \lambda ) = n \}$ form a partition of $\Lambda^\infty$ 
(see \ref{rad}(v)); moreover for every $\lambda \in \Lambda$ we have
\begin{equation} \label{decompose}
Z ( \lambda ) = \bigcup_{\substack{r ( \mu ) = s ( \lambda ) \\ d ( \mu ) = n}} 
Z ( \lambda \mu ) .
\end{equation}
\noindent
We endow $\Lambda^\infty$ with the topology generated by the collection 
$\{ Z (\lambda ) : \lambda \in \Lambda \}$. 
Note that the map given by $\lambda x \mapsto x$ induces a
homeomorphism between $Z ( \lambda )$ and $Z ( s ( \lambda ))$ for
all $\lambda \in \Lambda$. Hence, for every $p \in {\bf N}^k$ the map
$\sigma^p : \Lambda^\infty \rightarrow \Lambda^\infty$ is a local
homeomorphism.
\end{rmks}

\begin{lem} \label{zlambdacompact}
For each $\lambda \in \Lambda$, $Z ( \lambda )$ is compact.
\end{lem}

\begin{proof}
By \ref{partition} it suffices to show that $Z (v)$ is compact for all 
$v \in \Lambda^0$.  Fix $v \in \Lambda^0$ and let 
$\{ x_n \}_{n \ge 1}$ be a sequence in $Z ( v )$.   For every $m$, 
$x_n ( 0 , m )$ may take only finitely many values (by
(\ref{stand})).  Hence there is a $\lambda \in \Lambda^m$ such that 
$x_n( 0 , m ) = \lambda$ for infinitely many $n$. We may therefore
inductively construct a sequence $\{\lambda_j  : j \ge 1 \}$ in 
$\Lambda^p$ such that $s ( \lambda_j ) = r ( \lambda_{j+1} )$ and 
$x_n  ( 0 , jp  ) = \lambda_1 \cdots \lambda_j$ for infinitely many
$n$ (recall $p = (1 , \ldots , 1 ) \in {\bf N}^k$). 
Choose a  subsequence $\{ x_{n_j} \}$ such that
$x_{n_j} ( 0 , jp ) = \lambda_1 \cdots \lambda_j$.   Since $\{ jp \}$
is cofinal,  there is a unique $y \in \Lambda^\infty (v)$ such that 
$y ( (j-1)p , jp ) = \lambda_j$ for $j \ge 1$; then 
$x_{n_j} \rightarrow y$ and hence $Z(v)$ is compact.
\end{proof}

\noindent
Note that $\Lambda^\infty$ is compact if and only if $\Lambda^0$ is finite.

\begin{dfn} \label{gpdef}
If $\Lambda$ is $k$-graph then let 
\[
{\mathcal G}_{\Lambda} =  \{ ( x , n , y ) \in 
\Lambda^\infty \times {\bf Z}^k \times \Lambda^\infty : 
\sigma^\ell x = \sigma^m y , n = \ell - m\}. 
\]
\noindent
Define range and source maps  $r, s : {\mathcal G}_\Lambda \rightarrow
\Lambda^\infty$ by $r (x , n , y ) = x$, $s ( x , n , y ) = y$.
For  $( x, n, y )$, $( y, \ell , z ) \in {\mathcal G}_\Lambda$ set 
$( x , n , y ) (y, \ell , z ) = ( x , n + \ell , z )$,  and 
$( x, n, y )^{-1} =  ( y , -n , x )$; ${\mathcal G}_\Lambda$ is called
the path groupoid of $\Lambda$ (cf.~\cite{r,d,kprr}).
\end{dfn}

\noindent
One may check that ${\mathcal G}_\Lambda$ is a groupoid with  
$\Lambda^\infty = {\mathcal G}_\Lambda^0$ under the identification
$x \mapsto ( x, 0, x )$. For $\lambda $, $\mu \in \Lambda$
such that $s ( \lambda ) = s ( \mu )$ define 
\[
Z ( \lambda , \mu ) =  \{ (\lambda z,  d ( \lambda ) - d ( \mu ), \mu z ) : 
                                         z \in \Lambda^\infty(s (\lambda ))\}.
\]

\noindent
We collect certain standard facts about ${\mathcal G}_\Lambda$
in the following result:

\begin{prop} 
Let $\Lambda$ be a $k$--graph. The sets $\{ Z ( \lambda , \mu ) :
\lambda , \mu \in \Lambda , s ( \lambda ) = s ( \mu ) \}$ form a basis
for a locally compact Hausdorff topology on ${\mathcal G}_\Lambda$. With
this topology ${\mathcal G}_\Lambda$ is a second countable,
$r$--discrete locally compact groupoid in which each $Z ( \lambda
, \mu )$ is a compact open bisection. The topology on $\Lambda^\infty$ agrees
with the relative topology under the identification of $\Lambda^\infty$ with the subset
${\mathcal G}_\Lambda^0$ of ${\mathcal G}_\Lambda$. 
\end{prop}

\begin{proof}
One may check that the sets $Z ( \lambda , \mu )$ form a  basis for
a topology on ${\mathcal G}_\Lambda$. To see that multiplication is
continuous, suppose that $(x , n , y ) ( y , \ell , z ) = ( x , n + \ell , z )
\in Z ( \gamma ,\delta )$. Since $( x , n , y ), ( y , \ell , z )$ are
composable in ${\mathcal G}_\Lambda$ there are $\kappa , \nu \in
\Lambda$ and $t \in \Lambda^\infty$ such that $x = \gamma \kappa
t$, $y = \nu t$ and $z = \delta \kappa t$. Hence $( x , k , y ) \in Z (
\gamma \kappa , \nu )$ and $( y , \ell , z ) \in Z ( \nu , \delta
\kappa )$ and the product maps the open set ${\mathcal G}_\Lambda^2
\cap ( Z ( \gamma \kappa , \nu ) \times Z ( \nu , \delta \kappa ) )$
into $Z ( \gamma , \delta )$. The remaining parts of the proof are
similar to those given in \cite[Proposition 2.6]{kprr}.
\end{proof}
\noindent
Note that $Z ( \lambda , \mu ) \cong Z (s ( \lambda ) )$, via the map
$(\lambda z,  d ( \lambda ) - d ( \mu ), \mu z )  \mapsto z$.
Again we note that in the case $k=1$ we have $\Lambda = 
E^*$ for some directed graph $E$ and the groupoid
${\mathcal G}_{E^*} \cong {\mathcal G}_{E}$, the graph groupoid
of $E$ which is described in detail in \cite[\S 2]{kprr}.  

\begin{prop} \label{fstarprops}
Let $\Lambda$ be a $k$-graph and let $f : {\bf N}^\ell \rightarrow
{\bf N}^k$ be a morphism.  The map $x \mapsto f^* (x)$ given by  
$f^* (x) (m,n) = ( x (  f(m) , f (n) ) , n-m )$ defines a  continuous
surjective map 
$f^* : \Lambda^\infty \rightarrow f^* ( \Lambda )^\infty$. 
Moreover, if the image  of $f$ is cofinal (equivalently $f (p)$ is
strictly positive in the sense that all of its coordinates are
nonzero) then $f^*$ is a homeomorphism. 
\end{prop}

\begin{proof}
Given $x \in f^* ( \Lambda )^\infty$ choose a sequence $\{ m_i \}$
such that  $n_j = \sum_{i=1}^j m_i$ is cofinal in ${\bf N}^\ell$. 
Set $n_0 =0$ and let $\lambda_j \in \Lambda^{f(m_j)}$ be
defined by the condition that  
$x (n_{j-1} , n_j ) = ( \lambda_j, m_j )$.
We must show that there is an $x' \in \Lambda^\infty$ such that 
$x' (f (n_{j-1}) , f(n_j ) ) = \lambda_j$.  
It suffices to show that the the intersection 
$\cap_j Z ( \lambda_1 \cdots \lambda_j ) \neq \emptyset$. But this
follows by the finite intersection property.  One checks that 
$x = f^* (x')$. Furthermore the inverse image of $Z( \lambda, n )$ is
$Z ( \lambda )$ and hence $f^*$ is continuous.

Now suppose that the image  of $f$ is cofinal, then the procedure
defined above gives a continuous inverse for $f^*$.  Given 
$x \in f^* ( \Lambda )^\infty$, then since $f (n_j )$ is cofinal, the
intersection $\cap_j Z ( \lambda_1 \cdots \lambda_j )$ contains a
single point $x'$.  Note that $x'$ depends on $x$ continuously.
\end{proof}

\noindent
For higher rank graphs of the form $f^*(\Lambda )$ with $f$
surjective (see \ref{pullback}), the associated groupoid 
${\mathcal G}_{f^*(\Lambda )}$ decomposes as a direct product as
follows: 

\begin{prop} \label{groupoidpull}
Let $\Lambda$ be a $k$-graph and let $f : {\bf N}^\ell \rightarrow
{\bf N}^k$ be a surjective morphism. Then 
\[
 {\mathcal G}_{f^* ( \Lambda )} \cong
{\mathcal G}_{\Lambda} \times {\bf Z}^{\ell - k}   .
\]
\end{prop}

\begin{proof}
Since $f$ is surjective, the map 
$f^* : \Lambda^\infty \rightarrow f^* ( \Lambda )^\infty$  
is a homeomorphism (see \ref{fstarprops}).
The map $f$ extends to a surjective morphism
$f : {\bf Z}^\ell \rightarrow {\bf Z}^k$.  Let  
$j : {\bf Z}^k \rightarrow {\bf Z}^\ell$ be a section for
$f$ and let
$i : {\bf Z}^{\ell -k} \rightarrow {\bf Z}^\ell$ be an identification of
${\bf Z}^{\ell -k}$ with $\mbox{ker} \, f$.
Then we get a groupoid isomorphism by the map
\[
( (x,n,y) , m ) \mapsto ( f^* x , i (m) + j (n) , f^* y )  ,
\]
where $( (x,n,y), m ) \in {\mathcal G}_{\Lambda} \times {\bf Z}^{\ell -k}$.
\end{proof}

\noindent
Finally, as in \cite[Lemma 3.8]{rs2} we demonstrate that there is a nontrivial
$*$--representation of $( \Lambda , d)$.

\begin{prop} \label{univexists}
Let $\Lambda$ be $k$-graph  then there
exists a representation $\{ S_\lambda : \lambda \in \Lambda \}$ of
$\Lambda$ on a  Hilbert space with all partial isometries
$S_\lambda$ nonzero.
\end{prop}

\begin{proof}
Let ${\mathcal H} = \ell^2 ( \Lambda^\infty )$, then for $\lambda \in
\Lambda$ define $S_\lambda \in {\mathcal B} ( {\mathcal H} )$ by
\[
S_\lambda e_y ~=~ \left\{
\begin{array}{ll}
e_{\lambda y} ~&~ \mbox{if} ~ s ( \lambda ) = y (0) , \\
0 ~&~ \mbox{otherwise,}
\end{array}
\right.
\]
\noindent
where $\{ e_y : y \in \Lambda^\infty \}$ is the canonical basis for
${\mathcal H}$. Notice that $S_\lambda$ is nonzero since
$\Lambda^\infty ( s ( \lambda ) ) \neq \emptyset$;  one then checks
that the family $\{ S_\lambda : \lambda \in \Lambda \}$ satisfies
conditions \ref{cstarlambdadef}(i)--(iv). 
\end{proof}

\section{The gauge invariant uniqueness theorem}

\noindent
By the universal property of $C^* ( \Lambda )$  there is a
canonical action of the $k$-torus ${\bf T}^k$, called the {\bf gauge
action}: $\alpha : {\bf T}^k
\rightarrow \mbox{Aut}\, C^* ( \Lambda )$ defined for
$t = ( t_1 , \ldots , t_k ) \in {\bf T}^k$ and $s_\lambda \in C^*
( \Lambda )$ by
\begin{equation} \label{gauge}
\alpha_t ( s_\lambda ) = t^{d ( \lambda )} s_\lambda 
\end{equation}
\noindent
where $t^{m} =  t_1^{m_1} \cdots t_k^{m_k}$ for 
$m = ( m_1 , \ldots , m_k ) \in {\bf N}^k$.  
It is straightforward to show that $\alpha$ is strongly continuous.   
As in \cite[Lemma 2.2]{ck} and \cite[Lemma 3.6]{rs2} we shall need 
the following:

\begin{lem} \label{basicstructure}
Let $\Lambda$ be a $k$-graph. Then for $\lambda , \mu \in \Lambda$ and 
$q \in {\bf N}^k$ with $d ( \lambda )$, $d ( \mu ) \le q$ we have
\begin{equation} \label{starprod}
s_\lambda^* s_\mu ~=~
\sum_{\substack{\lambda\alpha = \mu\beta \\ d ( \lambda\alpha ) = q}}
s_\alpha  s_\beta^*.
\end{equation}
\noindent
Hence every nonzero word in $s_ \lambda , s_\mu^*$ may
be written as a finite sum of partial isometries of the form
$s_\alpha s_\beta^*$ where $s ( \alpha ) =  s ( \beta )$; their linear
span then forms a dense $*$--subalgebra of $C^* ( \Lambda )$.
\end{lem}

\begin{proof}
Applying \ref{cstarlambdadef}(iv) to $s ( \lambda )$ with
$n =  q - d ( \lambda )$, to $s ( \mu )$ with $n = q - d ( \mu )$
and using \ref{cstarlambdadef} (ii) we get
\begin{eqnarray}
s_\lambda^* s_\mu &=& p_{s ( \lambda )} s_\lambda^* s_\mu p_{s( \mu)} 
= \left(\sum_{
\Lambda^{q - d ( \lambda )} ( s ( \lambda ) )} \hspace{-5mm}
s_\alpha s_\alpha^* \right) s_\lambda^* s_\mu 
\left( \sum_{
\Lambda^{q - d ( \mu )} ( s ( \mu ) )} \hspace{-5mm} s_\beta
s_\beta^*  \right)
\nonumber \\
&=& \left( \sum_{
\Lambda^{q - d ( \lambda )} ( s ( \lambda ) )} \hspace{-5mm}
s_\alpha s_{\lambda \alpha}^* \right)
\left( \sum_{
\Lambda^{q - d ( \mu )} ( s ( \mu ) )}
\hspace{-5mm} s_{\mu \beta} s_\beta^* \right) .
 \label{lmdecompose}
\end{eqnarray}
\noindent
By \ref{rad}(iv) if $d ( \lambda \alpha ) = d ( \mu \beta )$ but
$\lambda \alpha \neq \mu \beta$, then the range projections
$p_{\lambda \alpha}$,  $p_{\mu \beta}$ are orthogonal and hence
one has $s_{\lambda \alpha}^* s_{\mu \beta}=  0$.  
If $\lambda \alpha = \mu \beta$ then 
$s_{\lambda \alpha}^* s_{\mu \beta} =  p_v$ where $v = s(\alpha)$
and so $s_\alpha s_{\lambda \alpha}^*s_{\mu \beta} s_\beta^* =
s_\alpha p_vs_\beta^* = s_\alpha s_\beta^*$; 
formula (\ref{starprod}) then follows from formula (\ref{lmdecompose}). 
%by applying \ref{cstarlambdadef}(i) to (\ref{lmdecompose}). 
The rest of the proof is now routine. 
%and we refer the reader to \cite[Lemma 1.1]{kpr} for the details.
\end{proof}

\noindent
Following \cite[\S 4]{rs2}: for $m \in {\bf N}^k$ let ${\mathcal F}_m$
denote the $C^*$--subalgebra of $C^* ( \Lambda )$ generated by the
elements $s_\lambda s_\mu^*$ for $\lambda , \mu \in
\Lambda^m$  where $s ( \lambda ) = s ( \mu )$, and for $v \in
\Lambda^0$ denote ${\mathcal F}_m (v)$ the $C^*$--subalgebra
generated by $s_\lambda s_\mu^*$  where $s ( \lambda ) = v$.

\begin{lem} \label{buildingblocks}
For $m \in {\bf N}^k$, $v \in \Lambda^0$ there exist 
isomorphisms 
$${\mathcal F}_m (v) \cong {\mathcal K} \left( \ell^2 (
\{ \lambda \in \Lambda^m : s ( \lambda ) = v \} ) \right)$$ 
and 
${\mathcal F}_m ~\cong~ \bigoplus_{v \in \Lambda^0} {\mathcal F}_m (v)$.
Moreover, the $C^*$--algebras ${\mathcal F}_m,$ $m \in {\bf N}^k$, 
form a directed system under inclusion, and
${\mathcal F}_\Lambda =\overline{ \cup {\mathcal F}_m}$ is an AF 
$C^*$--algebra.
\end{lem}

\begin{proof}
Fix $v \in \Lambda^0$ and let $\lambda$, $\mu$, $\alpha$, $\beta \in
\Lambda^m$ be such that $s ( \lambda ) = s ( \mu )$ and $s (\alpha )
= s ( \beta )$, then by \ref{rad}(iv) we have
\begin{equation} \label{mxunits}
\left( s_\lambda s_\mu^* \right) \left( s_\alpha s_\beta^*
\right) ~=~ \delta_{\mu , \alpha} s_\lambda s_\beta^* ,
\end{equation}
\noindent
so that the map which sends $s_\lambda s_\mu^* \in {\mathcal F}_m(v)$ 
to the matrix unit $e^v_{\lambda , \mu} \in {\mathcal K} 
\left( \ell^2 ( \{ \lambda \in \Lambda^m : s(\lambda ) = v \})\right)$ 
for all $\lambda$, $\mu \in \Lambda^m$ with $s(\lambda) = s(\mu) = v$
extends to an isomorphism. The second isomorphism also follows from
(\ref{mxunits}) (since $s(\mu) \neq s(\alpha)$ implies $\mu \neq \alpha$).
%from \ref{rad}(v) because for 
%$v , w \in \Lambda^0$ and $\lambda$, $\mu \in \Lambda^m$ are such that 
%$s( \lambda) = s (\mu ) = v$ and, $\alpha$, $\beta \in \Lambda^m$ are such that 
%$s( \alpha ) = s (\beta ) = w$ we have 
%$( s_\lambda s_\mu^* ) (s_\alpha s_\beta^* ) = 0$ if $v \neq w$. 
We claim that ${\mathcal F}_m$ is contained in ${\mathcal F}_{n}$ 
whenever $m \leq n$. To see this we apply \ref{cstarlambdadef}(iv) to give
\begin{equation} \label{embed}
s_\lambda s_\mu^* = s_\lambda p_{s ( \lambda )} s_\mu^* = 
\sum_{\Lambda^{\ell}(s (\lambda))} 
s_{\lambda} s_{\gamma} s_{\gamma}^*s_{\mu}^* =
\sum_{\Lambda^{\ell}(s (\lambda))} s_{\lambda \gamma} s_{\mu \gamma}^*
\end{equation}
\noindent
where $\ell = n - m$.  Hence the $C^*$--algebras 
${\mathcal F}_m,$ $m \in {\bf N}^k$, 
form a directed system as required. 
\end{proof}

\noindent
Note that ${\mathcal F}_\Lambda$ may also be expressed as the
closure of $\cup_{j=1}^\infty {\mathcal F}_{jp}$ where 
$p = ( 1 , \ldots , 1 ) \in {\bf N}^k$.

Clearly for $t \in {\bf T}^k$ the gauge automorphism $\alpha_t$
defined in (\ref{gauge}) fixes those elements 
$s_\lambda s_\mu^* \in C^* ( \Lambda )$
with $d ( \lambda ) = d ( \mu)$  (since $\alpha_t(s_\lambda s_\mu^*) = 
t^{d (\lambda ) - d ( \mu)}s_\lambda s_\mu^*$) 
and hence ${\mathcal F}_\Lambda$ is contained in the fixed point
algebra $C^* ( \Lambda )^\alpha$.  
Consider the linear map on $C^* (\Lambda )$ defined by
\[
\Phi (x) ~=~ \int_{{\bf T}^k} \alpha_t (x) \, dt
\]
\noindent
where $dt$ denotes normalised Haar measure on ${\bf T}^k$ and note
that $\Phi (x) \in C^* ( \Lambda )^\alpha$ for all 
$x \in C^* ( \Lambda )$.
As the proof of the following result is now standard, we omit it (see 
\cite[Proposition 2.11]{ck}, \cite[Lemma 3.3]{rs2}, \cite[Lemma
2.2]{bprs}):

\begin{lem}
Let $\Phi$, ${\mathcal F}_\Lambda$ be as described above.
\begin{itemize}
\item[(i)] The map $\Phi$ is a faithful conditional expectation 
from $C^* ( \Lambda )$ onto $C^* ( \Lambda )^\alpha$.
\item[(ii)] ${\mathcal F}_\Lambda = C^* ( \Lambda )^\alpha$.
\end{itemize}
\end{lem}

\noindent
Hence the fixed point algebra $C^* ( \Lambda )^\alpha$ is an AF 
algebra.
%Recall that by the universal property of $C^* ( \Lambda )$ given in
%\ref{universal} if $B$ is any $C^*$--algebra generated by a
%representation $\{ S_\lambda : \lambda \in \Lambda \}$ of $\Lambda$,
%then there is a $*$--homomorphism $\pi : C^* ( \Lambda ) \rightarrow
%B$ such that $\pi ( s_\lambda ) = S_\lambda$ for all $\lambda \in
%\Lambda$.
This fact is key to the proof of the gauge--invariant uniqueness  theorem
for $C^* ( \Lambda )$ (see \cite[Theorem 2.1]{bprs}, \cite[Theorem
2.3]{ahr}, see also \cite{ck,rs2} where a similar technique is used in
the proof of simplicity).

\begin{thm} \label{giuniqueness}
Let $B$ be a $C^*$--algebra, $\pi : C^* ( \Lambda ) \rightarrow B$
be a homomorphism and let 
$\beta : {\bf T}^k \rightarrow \mbox{Aut}\, (B)$ be an action such that 
$\pi\circ\alpha_t = \beta_t \circ \pi$ for all $t \in {\bf T}^k$.  Then $\pi$ 
is faithful if and only if $\pi ( p_v ) \neq 0$ for all $v \in \Lambda^0$.
\end{thm}

\begin{proof}
If $\pi ( p_v ) = 0$ for some $v \in\Lambda^0$ then clearly $\pi$
is not faithful.  Conversely, suppose that $\pi$ is equivariant and that 
$\pi ( p_v)\neq 0$ for all $v \in \Lambda^0$; we first  show that
$\pi$ is faithful on  
$C^*(\Lambda)^\alpha = \overline{\bigcup_{j \ge 0}
{\mathcal F}_{jp}}$. 
For any ideal $I$ in $C^* ( \Lambda )^\alpha$, we have 
$I = \overline{\bigcup_{j \ge 0} ( I \cap {\mathcal F}_{jp} )}$ 
(see \cite[Lemma 3.1]{b},
\cite[Lemma 1.3]{alnr}). Thus it is enough to prove that
$\pi$ is  faithful on each ${\mathcal F}_{n}$. But by \ref{buildingblocks}
it suffices to show that it is faithful on  
${\mathcal F}_{n} (v)$, for all $v \in \Lambda^0$.  Fix
$v \in \Lambda^0$ and $\lambda , \mu \in \Lambda^n$ with
$s ( \lambda ) = s ( \mu ) = v$ we need only show that 
$\pi ( s_\lambda s_\mu^* ) \neq 0$. Since $\pi ( p_v ) \neq 0$ we 
have
\[
0 \neq \pi ( p_v ^2 ) = \pi ( s_\lambda^* s_\lambda s_\mu^* s_\mu )
= \pi ( s_\lambda^* ) \pi ( s_\lambda s_\mu^* ) \pi ( s_\mu ) .
\]
Hence $\pi ( s_\lambda s_\mu^* ) \neq 0$ and $\pi$ is faithful on 
$C^* ( \Lambda)^\alpha$.  Let $a \in C^* ( \Lambda)$ be a  nonzero
positive element; then since $\Phi$ is faithful $\Phi (a) \neq 0$ 
and as $\pi$ is faithful on $C^* (\Lambda^\alpha )$ we have
$$
0 \neq \pi(\Phi (a)) = \pi\left(\int_{{\bf T}^k} \alpha_t (a) \, dt \right)
                                   = \int_{{\bf T}^k} \beta_t (\pi(a)) \, dt;
$$
hence, $\pi(a) \neq 0$ and $\pi$ is faithful on $C^* ( \Lambda)$ as required.
%By \cite[Lemma 2.2]{bkr} the faithfulness of $\pi$ will now
%follow provided that we show 
%$\| \pi ( \Phi (a) ) \| \le \| \pi (a) \|$ for all $a \in C^* ( \Lambda )$;
%fix $a \in C^* ( \Lambda )$ then by averaging over the strongly 
%continuous action $\beta$ we have
%\[
%\| \pi ( \Phi (a) ) \| \le \int_{{\bf T}^k} \| 
%\pi ( \alpha_t (a) ) \| dt = \int_{{\bf T}^k} \|
%\beta_t ( \pi (a) ) \| dt = \int_{{\bf T}^k}  \| \pi (a)
%\| dt =  \| \pi (a) \| ,
%\]
%\noindent
%as required. 
\end{proof}

\begin{cor} \label{giutconsequences} 
\hskip1ex
\begin{enumerate}
\item[(i)]
Let $( \Lambda , d )$ be a $k$-graph and let ${\mathcal G}_\Lambda$ be
its associated groupoid, then there is an isomorphism 
$C^* ( \Lambda ) \cong C^* ( {\mathcal G}_\Lambda )$ such that 
$s_\lambda\mapsto 1_{Z ( \lambda , s ( \lambda ) )}$ for $\lambda \in\Lambda$.
Moreover the canonical map 
$C^* ( {\mathcal G}_\Lambda ) \rightarrow C^*_r ({\mathcal G}_\Lambda )$ is an
isomorphism.
\item[(ii)]
Let $\{ M_1 , \ldots , M_k \}$ be a collection of
matrices satisfying (H0)--(H3) of \cite{rs2} and $W$ the $k$-graph defined in
\ref{secondex}(iv), then
$C^* ( W ) \cong {\mathcal A}$, via the map 
$s_\lambda \mapsto s_{\lambda , s ( \lambda )}$ for $\lambda \in W$.
\item[(iii)] If $\Lambda$ is a $k$-graph and $f : {\bf N}^\ell
\rightarrow {\bf N}^k$ is injective then the $*$-homomorphism $\pi_f :
C^* ( f^* ( \Lambda ) ) \rightarrow C^* ( \Lambda )$ (see 
\ref{pifprop}) is injective. In particular the $C^*$--algebras of the
coordinate graphs $\Lambda_i$ for $1 \le i \le k$ form a generating
family of subalgebras of $C^* ( \Lambda )$.
Moreover, if $f$ is surjective then 
$C^* (f^* ( \Lambda ) ) \cong 
C^* ( \Lambda ) \otimes C ( {\bf T}^{\ell - k})$.
\item[(iv)]
Let $( \Lambda_i , d_i )$ be $k_i$-graphs for 
$i = 1, 2$, then $C^* ( \Lambda_1 \times \Lambda_2 ) \cong
C^* ( \Lambda_1 ) \otimes C^* ( \Lambda_2 )$ via the map
$s_{(\lambda_1 , \lambda_2)} \mapsto s_{\lambda_1} \otimes
s_{\lambda_2}$ for $( \lambda_1 , \lambda_2 ) \in \Lambda_1
\times \Lambda_2$.
\end{enumerate}
\end{cor}

\begin{proof}
For (i) we note that 
$s_\lambda\mapsto 1_{Z ( \lambda , s ( \lambda ) )}$ 
for $\lambda \in \Lambda$ is a $*$-representation of $\Lambda$; hence
there is a $*$-homomorphism 
$\pi : C^* ( \Lambda ) \rightarrow C^* ( {\mathcal G}_\Lambda )$ 
such that
$\pi (s_\lambda ) = 1_{Z ( \lambda , s ( \lambda ))}$ 
for $\lambda \in \Lambda$ (see \ref{rad}(i)).
Let $\beta$ denote the ${\bf T}^k$-action on 
$C^* ( {\mathcal G}_\Lambda )$ induced by the ${\bf Z}^k$-valued 
$1$--cocycle defined on ${\mathcal G}_\Lambda$ by 
$(x,k,y) \mapsto k$ (see \cite[II.5.1]{r}); one checks that   
$\pi\circ\alpha_t = \beta_t \circ \pi$ for all $t \in {\bf T}^k$. 
Clearly for $v \in \Lambda^0$ we have $1_{Z(v,v)} \neq 0$,
since $\Lambda^\infty (v) \neq \emptyset$ and $\pi$ is injective.
Surjectivity follows from the fact that 
$\pi ( s_\lambda s_\mu^* ) = 1_{Z (\lambda , \mu )}$  together with 
the observation that $C^* ( {\mathcal G}_\Lambda ) = 
\overline{\mbox{span}}\{ 1_{Z(\lambda , \mu )} \}$. 
The same argument shows that $C^*_r ( {\mathcal G}_\Lambda )
\cong C^* ( \Lambda )$ and so $C^*_r ( {\mathcal G}_\Lambda ) \cong  C^* (
{\mathcal G}_\Lambda )$\footnote{This can be also deduced from the amenability
of ${\mathcal G}_\Lambda$ (see \ref{Gamenable})}.

For (ii) we note that there is a surjective $*$-homomorphism  
$\pi : C^*( W ) \rightarrow {\mathcal A}$ such that  
$\pi ( s_\lambda ) = s_{\lambda , s (\lambda )}$ for $\lambda \in W$
(see \ref{secondex}(iv))
which is clearly equivariant for the respective ${\bf T}^k$--actions. 
Moreover by \cite[Lemma 2.9]{rs2} we have
$s_{v,v} \neq 0$ for all $v \in W_0 = A$ and so the result follows

For (iii) note that the injection $f : {\bf N}^\ell \rightarrow {\bf N}^k$
extends naturally to a homomorphism $f : {\bf Z}^\ell \rightarrow {\bf
Z}^k$ which in turn induces a map $\hat{f} : {\bf T}^k \rightarrow
{\bf T}^\ell$ characterised by $\hat{f} (t)^p = t^{f(p)}$ for $p \in
{\bf N}^\ell$. Let $B$ be the fixed point algebra of the gauge
action of ${\bf T}^k$ on $C^* ( \Lambda )$ restricted to the kernel of
$\hat{f}$. The gauge action restricted to $B$ descends to an action of 
${\bf T}^\ell = {\bf T}^k / \mbox{Ker} \, \hat{f}$ on $B$ which we denote
$\overline{\alpha}$. Observe that for $t \in {\bf T}^k$ and
$(\lambda , n) \in  f^*(\Lambda)$ we have
\[
\alpha_t ( \pi_f ( s_{\lambda , n})) = t^{f(n)} s_\lambda =
\hat{f} (t)^n s_\lambda ;
\]
hence $\mbox{Im} \, \pi_f \subseteq B$ (if $t \in \mbox{Ker} \,
\hat{f}$ then $\hat{f} (t)^n =1$). By the same formula we see that
$\pi_f \circ\alpha = \overline{\alpha} \circ \pi_f$ and the result now
follows by \ref{giuniqueness}.  The last assertion follows from part (i) 
together with the fact that 
${\mathcal G}_{f^*(\Lambda)} \cong {\mathcal G}_\Lambda\times{\bf Z}^{\ell -k}$
(see \ref{groupoidpull}).

For (iv) define a map $\pi : C^* ( \Lambda_1 \times \Lambda_2 )
\rightarrow C^* ( \Lambda_1 ) \otimes C^* ( \Lambda_2 )$ 
given by $s_{( \lambda_1 , \lambda_2 )} \mapsto s_{\lambda_1}
\otimes s_{\lambda_2}$; this is surjective as these elements
generate $C^* ( \Lambda_1 ) \otimes C^* ( \Lambda_2 )$. 
We note that 
$C^* ( \Lambda_1 ) \otimes C^* (  \Lambda_2 )$ carries a 
${\bf T}^{k_1 +k_2}$ action $\beta$  defined for $( t_1 , t_2 ) \in 
{\bf T}^{k_1 + k_2}$ and 
$( \lambda_0 , \lambda_1 ) \in \Lambda_1 \times \Lambda_2$  by
$\beta_{( t_1 , t_2 )} ( s_{\lambda_1} \otimes s_{\lambda_2} ) =
\alpha_{t_1} s_{\lambda_1} \otimes \alpha_{t_2} s_{\lambda_2}$. 
Injectivity then follows by \ref{giuniqueness}, since $\pi$ is equivariant  
and for $( v , w )\in ( \Lambda_1 \times \Lambda_2 )^0$ we have  
$p_v \otimes p_w \neq 0$.
\end{proof}
\noindent
Henceforth we shall tacitly identify $C^* (\Lambda )$ with 
$C^* ( {\mathcal G}_\Lambda )$.
\begin{rmk}
Let $\Lambda$ be a $k$-graph and suppose that 
$f : {\bf N}^\ell \rightarrow {\bf N}^k$ is an injective morphism for
which $H$, the image of $f$, is cofinal.  Then $\pi_f$ induces an
isomorphism of $C^* ( f^* ( \Lambda ) )$ with its range, the fixed
point algebra of the restriction of the gauge action to $H^\perp$. 
\end{rmk}

\section{Aperiodicity and its consequences}

The aperiodicity condition we study in this section is an analog
of condition (L) used in \cite{kpr}. We first define
what it means for an infinite path to be periodic or aperiodic.

\begin{dfns}
For $x \in \Lambda^\infty$ and $p \in {\bf Z}^k$ we say that $p$
is {\bf a period} of $x$ if for every $(m,n) \in \Omega$ with $m+p
\ge 0$ we have $x ( m+p , n+p ) = x ( m , n )$. We say that $x$ is
{\bf periodic} if it has a nonzero period. We say that $x$ is
{\bf eventually periodic} if $\sigma^n x$ is periodic for some $n \in
{\bf N}^k$, otherwise $x$ is said to be {\bf aperiodic}.
\end{dfns}

\begin{rmks}
For $x \in \Lambda^\infty$ and $p \in {\bf Z}^k$, $p$
is {\bf a period} of $x$ if and only if $\sigma^m x = \sigma^n x$ for
all $m,n \in {\bf N}^k$ such that $p =m-n$. Similarly $x$ is
eventually periodic, with eventual period $p \neq 0$  if and only if 
$\sigma^m x = \sigma^n x$ for some $m,n \in {\bf N}^k$ such that
$p =m-n$.
\end{rmks}

\begin{dfn}
The $k$-graph $\Lambda$ is said to satisfy the {\bf aperiodicity condition} (A) if for
every $v \in \Lambda^0$ there is an aperiodic path $x \in\Lambda^\infty (v)$.
\end{dfn}

\begin{rmk}
Let $E$ be a directed graph which is row finite and has no sinks then
the associated $1$-graph $E^*$ satisfies the aperiodicity
condition if and only if every loop in $E$ has an exit (i.e.\ satisfies
condition (L) of \cite{kpr}). However, if we consider the $2$-graph
$f^* ( E^* )$ where $f : {\bf N}^2 \rightarrow {\bf N}$ is given
by $f ( m_1, m_2 ) = m_1 + m_2$ then $p=(1,-1)$ is a period
for every point in $f^* ( E^* )^\infty$ (even if $E$ has no loops).
\end{rmk}

\begin{prop} 
The groupoid ${\mathcal G}_\Lambda$ is essentially free (i.e. the points
with trivial isotropy are dense in ${\mathcal G}_\Lambda^0$) if and only
if $\Lambda$ satisfies the aperiodicity condition.
\end{prop}

\begin{proof}
Observe that if $x \in \Lambda^\infty$ is aperiodic then $\sigma^m x
= \sigma^n x$ implies that $m=n$ and hence 
$x \in \Lambda^\infty = {\mathcal G}_\Lambda^0$ has trivial isotropy, and
conversely. Hence ${\mathcal G}_\Lambda$ is essentially free if and only if
aperiodic points are dense in $\Lambda^\infty$. If aperiodic points are dense in
$\Lambda^\infty$ then $\Lambda$ clearly satisfies the aperiodicity condition, for
$Z(v) = \Lambda^\infty(v)$ must then contain aperiodic  points for every 
$v \in \Lambda^0$.  Conversely, suppose that $\Lambda$ satisfies the aperiodicity
condition, then for every $\lambda \in \Lambda$ there is 
$x \in \Lambda^\infty (s ( \lambda ))$ which is aperiodic. Then
$\lambda x  \in Z ( \lambda )$ is aperiodic. Hence the aperiodic points are dense in
$\Lambda^\infty$. 
\end{proof}

\noindent
The isotropy group of an element $x \in \Lambda^\infty$ is equal to
the subgroup of its eventual periods (including $0$).

\begin{thm}
Let $\pi : C^* ( \Lambda ) \rightarrow B$ be a $*$--homomorphism and suppose that
$\Lambda$ satisfies the aperiodicity condition. Then $\pi$ is faithful
if and only if 
$\pi ( p_v ) \neq 0$ for all $v \in \Lambda^0$.
\end{thm}

\begin{proof}
If $\pi ( p_v ) =0$ for some $v \in\Lambda^0$ then clearly $\pi$
is not faithful.  Conversely, suppose $\pi ( p_v ) \neq 0$ for all 
$v \in \Lambda^0$; then by \ref{giutconsequences}(i) we have
$C^* ( \Lambda ) = C^*_r ( {\mathcal G}_\Lambda )$ and hence from
\cite[Corollary 3.6]{kpr} it suffices to show that $\pi$ is 
faithful on $C_0 ( {\mathcal G}_\Lambda^0 )$. If the kernel  of
the restriction of $\pi$ to $C_0 ( {\mathcal G}_\Lambda^0 )$ is
nonzero, it must contain the characteristic function 
$1_{Z ( \lambda )}$ for some $\lambda \in \Lambda$. 
It follows that $\pi ( s_\lambda s_\lambda^* ) = 0$ and hence
$\pi ( s_\lambda ) = 0$; in which case 
$\pi ( p_{s ( \lambda )} ) =  \pi ( s_\lambda^* s_\lambda ) =0$,
a contradiction.
\end{proof}

\begin{dfn}
We say that $\Lambda$ is {\bf cofinal} if for every $x \in
\Lambda^\infty$ and $v \in \Lambda^0$ there is $\lambda \in
\Lambda$ and $n \in {\bf N}^k$ such that $s ( \lambda ) = x(n)$
and $r ( \lambda ) = v$.
\end{dfn}

\begin{prop} \label{simple}
Suppose $\Lambda$ satisfies the aperiodicity condition, then 
$C^*(\Lambda )$ is simple if and only if $\Lambda$ is cofinal.
\end{prop}

\begin{proof}
By \ref{giutconsequences}(i) 
$C^* ( \Lambda ) = C^*_r ( {\mathcal G}_\Lambda )$;
since $ {\mathcal G}_\Lambda$ is essentially free, 
$C^* ( \Lambda )$ is simple if and only if ${\mathcal G}_\Lambda$ 
is minimal. Suppose that $\Lambda$ is cofinal and fix
$x \in \Lambda^\infty$ and $\lambda \in \Lambda$; then by cofinality
there is a $\mu \in \Lambda$ and $n \in {\bf N}^k$ so that 
$s ( \mu ) = x(n)$ and $r ( \mu  ) =  s ( \lambda )$. Then
$y = \lambda \mu \sigma^n x \in Z ( \lambda )$ and $y$ is in the same
orbit as $x$; hence all orbits are  dense and ${\mathcal G}_\Lambda$
is minimal. 

Conversely, suppose that ${\mathcal G}_\Lambda$ is minimal and that 
$x \in \Lambda^\infty$ and $v \in \Lambda^0$ then there is 
$y \in Z (v)$ such that $x, y$ are in the same orbit. Hence there
exist $m, n \in {\bf N}^k$ such that $\sigma^n x = \sigma^m y$; 
then it is easy to check that $\lambda = y ( 0, m )$ and $n$ have 
the desired properties. 
\end{proof}

\noindent
Notice that second hypothesis used in the following corollary is the
analog of the condition that every vertex connects to a loop and it is
equivalent to requiring that for every $v \in \Lambda^0$, there is an
eventually periodic $x \in \Lambda^\infty (v)$ with positive eventual
period (i.e.\ the eventual period lies in ${\bf N}^k \backslash \{ 0 \}$).
The proof follows the same lines as \cite[Theorem 3.9]{kpr}:

\begin{prop} \label{pi}
Let $\Lambda$ satisfy the aperiodicity condition. Suppose that for
every $v \in \Lambda^0$ there are $\lambda$, $\mu \in \Lambda$ with 
$d(\mu) \neq 0$ such that $r ( \lambda ) = v$ and 
$s (\lambda ) = r ( \mu ) = s (\mu )$ then $C^*( \Lambda )$ is purely infinite 
in the sense that every hereditary subalgebra contains an infinite projection.
\end{prop}

\begin{proof}
Arguing as in \cite[Lemma 3.8]{kpr} one shows that ${\mathcal G}_\Lambda$
is locally contracting. The aperiodicity condition guarantees that 
${\mathcal G}_\Lambda$ is essentially free, hence by \cite[Proposition 2.4]{a-d}
(see also \cite{ls}) we have 
$C^* ( \Lambda ) = C_r^* ( {\mathcal G}_\Lambda )$ is purely infinite.
\end{proof}

\section{Skew products and group actions}

Let $G$ be a discrete group, $\Lambda$ a $k$-graph and  
$c :\Lambda \rightarrow G$ a functor.  We introduce an analog of the
skew product graph considered in \cite[\S 2]{kp} (see also \cite{gt});
the resulting object, which we denote $G \times_c \Lambda$, is also a
$k$-graph.   As in \cite{kp} if $G$ is abelian the associated 
$C^*$--algebra is isomorphic to a crossed product of $C^*(\Lambda)$ by
the natural action of $\widehat{G}$ induced by $c$ (more generally it
is a crossed product by a coaction --- see \cite{ma,kqr}). As a
corollary we show that the crossed product of $C^* (\Lambda )$ by the
gauge action, $C^*(\Lambda ) \rtimes_\alpha {\bf T}^k$, is isomorphic to  
$C^*({\bf Z}^k \times_d \Lambda)$, the $C^*$--algebra of the
skew-product $k$-graph arising from the degree map. It will then
follow that $C^* ( \Lambda) \rtimes_\alpha {\bf T}^k$ is AF and that 
${\mathcal G}_\Lambda$ is amenable.   

\begin{dfn} 
Let $G$ be a discrete group, $( \Lambda , d)$ a $k$-graph. Given $c :
\Lambda \rightarrow G$ a functor then define the {\bf skew product}
$G \times_c \Lambda$ as follows: the objects are identified with 
$G \times \Lambda^0$ and the morphisms are identified with 
$G \times \Lambda$ with the following structure maps
\[
s ( g , \lambda ) = ( g c ( \lambda ) , s ( \lambda ) )  \quad
\mbox{and} \quad r ( g , \lambda ) = ( g , r ( \lambda ) ) .
\]
If $s ( \lambda ) = r ( \mu )$ then $(g , \lambda )$ and $( g c (
\lambda ) , \mu)$ are composable in $G \times_c \Lambda$ and
\[
( g , \lambda ) ( g c ( \lambda) , \mu ) = ( g , \lambda \mu ) .
\]
The degree map is given by $d ( g , \lambda ) = d ( \lambda )$.
\end{dfn}

\noindent
One must check that $G \times_c \Lambda$ is a $k$-graph.  
%is indeed a small category and that $d$ satisfies the factorisation property; 
%so $G \times_c \Lambda$ 
If $k=1$ then any function $c : E^1 \rightarrow G$ extends to a unique functor 
$c : E^* \rightarrow G$ (as in \cite[\S 2]{kp}).   
The skew product graph $E(c)$ of 
\cite{kp} is related to our  skew product in a simple way:  
$G \times_c E^* = E(c)^*$.
A key example of this construction arises by regarding the degree map
$d$ as a functor with values in ${\bf Z}^k$. 

The functor $c$ induces a cocycle 
$\tilde{c} : {\mathcal G}_\Lambda \rightarrow G$
as follows: given $(x, \ell - m , y ) \in {\mathcal G}_\Lambda$ so that 
$\sigma^\ell x =  \sigma^m y$ then set 
\[
\tilde{c} (x , \ell  - m , y ) = c ( x ( 0 , \ell ) ) c ( y ( 0 , m ) )^{-1} .
\]
\noindent
As in \cite{kp} one checks that this is well-defined and that
$\tilde{c}$ is a (continuous) cocycle; regarding the degree map
$d$ as a functor with values in ${\bf Z}^k$, we have 
$\tilde{d} (x ,  n , y ) = n$ for $(x ,  n , y ) \in  {\mathcal G}_\Lambda$.
In the following we show that the skew product groupoid obtained from 
$\tilde{c}$ (as defined in \cite{r}) is the same as the path groupoid of 
the skew product (cf.~\cite[Theorem 2.4]{kp}):

\begin{thm} \label{skewfunctor}
Let $G$ be a discrete group, $\Lambda$ a $k$-graph and 
$c : \Lambda \rightarrow G$ a functor. Then 
${\mathcal G}_{G \times_c \Lambda} \cong {\mathcal G}_\Lambda (\tilde{c})$
where $\tilde{c} : {\mathcal G}_\Lambda \rightarrow G$ is defined as above.
\end{thm}

\begin{proof}
We first identify $G \times \Lambda^\infty$ with 
$(G \times_c \Lambda )^\infty$  as
follows: for $( g , x ) \in G \times \Lambda^\infty$ define 
$(g, x) :  \Omega \rightarrow G \times_c \Lambda$ by 
\[
(g , x) (m,n) = ( g c ( x (0,m ) ) , x (m,n) ) ;
\]
it is straightforward to check that this defines a
degree--preserving functor and thus an element of 
$(G \times_c \Lambda )^\infty$. Under this identification $\sigma^n
(g,x) = ( g c ( x  ( 0 , n ) ) , \sigma^n x )$ for all $n \in {\bf N}^k$,
$(g,x) \in ( G \times_c \Lambda )^\infty$. As in the proof of
\cite[Theorem 2.4]{kp} define a map 
$\phi : {\mathcal G}_\Lambda ( \tilde{c} ) \rightarrow 
{\mathcal G}_{G \times_c \Lambda}$ as follows: for $x , y \in \Lambda^\infty$ with
$\sigma^\ell x = \sigma^m y$ set 
$\phi ( [ x ,\ell - m , y ] , g ) = ( x' ,\ell - m , y' )$ where $x' = ( g , x )$ and 
$y' = ( g \tilde{c} ( x , \ell - m , y ) , y )$. Note that
\begin{eqnarray*}
\sigma^m y' = \sigma^m ( g \tilde{c} ( x , \ell - m , y ) , y ) &=& 
 \sigma^m ( g c ( x ( 0 , \ell ) ) c ( y ( 0 , m ) )^{-1} , y ) \\
&=& ( g c ( x ( 0 , \ell ) ) , \sigma^m y ) = 
( g c ( x ( 0 , \ell ) ) , \sigma^\ell x ) = \sigma^\ell (g,x) = 
\sigma^\ell x' ,
\end{eqnarray*}
and hence $( x',\ell - m , y' ) \in {\mathcal G}_{G \times_c \Lambda}$.
The rest of the proof proceeds as in \cite[Theorem 2.4]{kp} {\em
mutatis mutandis}.
\end{proof}

\begin{cor} \label{after}
Let $G$ be a discrete abelian group, $\Lambda$ a $k$-graph and 
$c : \Lambda \rightarrow G$ a functor. There is an action 
$\alpha^c : \widehat{G} \rightarrow \hbox{\rm Aut} \, C^* ( \Lambda)$ 
such that for $\chi \in \widehat{G}$ and $\lambda \in \Lambda$
\[
\alpha^c_\chi  ( s_\lambda  ) = \langle \chi , c ( \lambda ) \rangle
s_\lambda .
\]
Moreover $C^* ( \Lambda ) \rtimes_{\alpha^c} \widehat{G} 
\cong C^*( G \times_c \Lambda )$. In particular the gauge action is of
the form, $\alpha = \alpha^d$, and so $C^* ( \Lambda ) \rtimes_\alpha
{\bf T}^k \cong C^* ( {\bf Z}^k \times_d \Lambda )$. 
\end{cor}

\begin{proof}
Since $C^* ( \Lambda )$ is defined to be the universal $C^*$--algebra
generated by the $s_\lambda$'s subject to the relations
(\ref{cstarlambdadef}) and $\alpha^c$ preserves these relations it is
clear that it defines an action of $\widehat{G}$ on $C^* ( \Lambda )$.
The rest of the proof follows in the same manner as that of
\cite[Corollary 2.5]{kp} (see \cite[II.5.7]{r}).
\end{proof}

\noindent
In order to show that $C^* ( \Lambda ) \rtimes_\alpha {\bf T}^k$ is
AF, we need the following lemma:

\begin{lem} \label{before}
Let $\Lambda$ be a $k$-graph and suppose there is a map 
$b : \Lambda^0 \rightarrow {\bf Z}^k$ such that 
$d ( \lambda ) = b ( s ( \lambda ) ) - b ( r ( \lambda ) )$ for all
$\lambda \in \Lambda$, then $C^* ( \Lambda )$ is AF. 
\end{lem}

\begin{proof}
For every $n \in {\bf Z}^k$ let $A_n$ be the closed linear span of
elements of the form $s_\lambda s_\mu^*$ with 
$b ( s ( \lambda ) ) = n$.  Fix $\lambda$, $\mu \in \Lambda$ with 
$b ( s ( \lambda ) ) = b ( s ( \mu ) ) = n$ we claim that
$s_\lambda^* s_\mu = 0$ if $\lambda \neq \mu$. 
If $s_\lambda^* s_\mu \neq 0$ then by  \ref{basicstructure} there are 
$\alpha$,  $\beta \in \Lambda$ with $s(\lambda ) = r(\alpha)$ and 
$s(\mu ) = r(\beta)$ such that $\lambda \alpha = \mu \beta$; but then we have
$$ 
d (\alpha ) + n  = d (\alpha ) + b ( s ( \lambda ) ) = 
b ( s ( \lambda \alpha ) ) = b ( s ( \mu \beta ) )  = 
d(\beta) + b (s ( \mu ) ) = d(\beta) + n.
$$ 
Thus $d ( \alpha ) = d ( \beta )$ and hence by the factorisation
property $\alpha = \beta$. Consequently $\lambda = \mu$ by
cancellation and the claim is established.   It follows that for each
$v$ with $b(v)=n$ the elements $s_\lambda s_\mu^*$ with 
$s(\lambda ) = s (\mu ) = v$ form a system of matrix units and two
systems associated to distinct $v$'s are orthogonal (see 
\ref{buildingblocks}).  Hence we have
\[
A_n \cong \bigoplus_{b(v)=n} {\mathcal K} \left( \ell^2 ( s^{-1} (v) \right).
\]
By an argument similar to that in the proof of Lemma \ref{buildingblocks},
if $n \leq m$ then $A_n \subseteq A_m$ (see equation (\ref{embed})); 
our conclusion now follows. 
\end{proof}

\noindent
Note that $A_n$ in the above proof is the $C^*$--algebra
of a subgroupoid of ${\mathcal G}_\Lambda$ which is isomorphic to the
disjoint union
$$
\bigsqcup_{b(v)=n} R_v \times \Lambda^\infty (v)
$$
where $R_v$ is the transitive principal groupoid on $s^{-1} (v)$. 
Since ${\mathcal G}_\Lambda$ is the increasing union of these
elementary groupoids, it is an AF-groupoid and hence amenable 
(see \cite[III.1.1]{r}). The existence of such a function 
$b  : \Lambda^0 \rightarrow {\bf Z}^k$ is not 
necessary for $C^* ( \Lambda )$ to be AF since there are $1$--graphs with no loops
which do not have this property (see \cite[Theorem 2.4]{kpr}).

\begin{thm} \label{Gamenable}
Let $\Lambda$ be a $k$-graph, then $C^* ( \Lambda ) \rtimes_\alpha
{\bf T}^k$ is AF and the groupoid ${\mathcal G}_\Lambda$ is amenable. 
Moreover, $C^* ( \Lambda )$ falls in the bootstrap class ${\mathcal N}$ of
\cite{rs-uct} and is therefore nuclear.  
Hence, if $C^*(\Lambda)$ is simple and purely infinite (see \S 5), then it
may be classified by its $K$-theory.
\end{thm}

\begin{proof}
Observe that the map $b : ( {\bf Z}^k \times_d \Lambda )^0
\rightarrow {\bf Z}^k$ given by $b ( n , v ) = n$ satisfies  
$$
b ( s ( n , \lambda ) ) -  b ( r ( n ,\lambda )) =
b(n +  d ( \lambda ), \lambda ) -  b ( n , r(\lambda )) = 
n +  d ( \lambda ) - n = d ( n ,\lambda )  . 
$$
The first part of the result then follows from
\ref{before} and \ref{after}. 
To show that ${\mathcal G}_\Lambda$ is amenable 
we first observe that 
${\mathcal G}_{\Lambda} ( \tilde{d} ) \cong 
{\mathcal G}_{{\bf Z}^k \times_d \Lambda}$ 
is amenable. Since ${\bf Z}^k$ is amenable,
% $\tilde{d} : {\mathcal G}_\Lambda \rightarrow {\bf Z}^k$ 
we may apply \cite[Proposition II.3.8]{r} to deduce that 
${\mathcal G}_\Lambda$ is amenable. % (see also \cite[p. 65]{a-dr}). 
Since $C^* ( \Lambda )$ is strongly Morita equivalent to the
crossed product of an AF algebra by a ${\bf Z}^k$--action, it falls in the 
bootstrap class ${\mathcal N}$ of \cite{rs-uct}.
The final assertion follows from the Kirchberg-Phillips classification theorem 
(see \cite{k,p}). 
\end{proof}

\noindent
We now consider free actions of groups on $k$-graphs (cf.\
\cite[\S 3]{kp}). Let $\Lambda$ be a $k$-graph and $G$ a countable
group, then $G$ {\bf acts on} $\Lambda$ if there is a group
homomorphism  $G \rightarrow\mbox{Aut} \,\Lambda$ (automorphisms
are compatible with all structure maps, including the degree): write 
$(g , \lambda ) \mapsto g  \lambda$. The action of $G$ on
$\Lambda$ is said to be {\bf free} if it is free on $\Lambda^0$.
By the universality of $C^* (\Lambda )$ an action of $G$ on
$\Lambda$ induces an action $\beta$ on $C^* ( \Lambda )$ such that
$\beta_g s_\lambda = s_{g \lambda}$.

Given a free action of a group $G$ on a $k$-graph $\Lambda$ one
forms the {\bf quotient} $\Lambda / G$ by the equivalence relation
$\lambda \sim \mu$ if $\lambda = g \nu$ for some $g \in G$. One
checks that all structure maps are compatible with $\sim$ and so
$\Lambda /G$ is also a $k$-graph. 

\begin{rmk} \label{gtanalog}
Let $G$ be a countable group and $c : \Lambda \rightarrow G$ a
functor, then $G$ acts freely on $G \times_c \Lambda$ by $g ( h ,
\lambda ) = ( gh , \lambda )$; furthermore $(G \times_c \Lambda ) / G
\cong \Lambda$. 

Suppose now that $G$ acts freely on $\Lambda$ with quotient 
$\Lambda / G$; we claim that $\Lambda$ is isomorphic, in an
equivariant way, to a skew product of $\Lambda / G$ for some 
suitably chosen $c$ (see \cite[Theorem 2.2.2]{gt}).  
Let $q$ denote the quotient map. For every $v \in ( \Lambda / G )^0$ 
choose $v' \in \Lambda^0$ with $q ( v' ) = v$ and for every 
$\lambda \in \Lambda / G$ let $\lambda'$ denote the unique element in 
$\Lambda$ such that $q ( \lambda' ) = \lambda$ and 
$r ( \lambda' ) = r ( \lambda)'$.  Now let $c :\Lambda / G \rightarrow G$ 
be defined by the formula
\[
s ( \lambda' ) =  c ( \lambda ) s ( \lambda )' .
\]
We claim that $c ( \lambda \mu ) = c ( \lambda ) c ( \mu )$
for all $\lambda$, $\mu \in \Lambda$ with $s(\lambda) = r(\mu)$.
Note that  
$$
r( c ( \lambda ) \mu' ) =  c ( \lambda ) r(\mu' ) = c ( \lambda ) r(\mu)' = 
                   c ( \lambda ) s( \lambda )' =  s( \lambda ');
$$
hence, we have
$( \lambda \mu)' = \lambda' ( c ( \lambda ) \mu' )$ 
(since the image of both sides
agree under $q$ and $r$).  Thus
\[
c (\lambda \mu ) s ( \mu )' = c ( \lambda \mu ) s ( \lambda \mu )' = 
s [ ( \lambda \mu )' ] = s (  c ( \lambda ) \mu' ) = 
c ( \lambda ) s (\mu' ) = c (\lambda ) c ( \mu ) s ( \mu )'
\]
\noindent
which establishes the desired identity (since $G$ acts freely on
$\Lambda$). The map $(g , \lambda )\mapsto g \lambda'$ defines an
equivariant isomorphism between $G \times_c ( \Lambda / G )$ and 
$\Lambda$ as required.
\end{rmk}

\noindent
The following is a generalization of \cite[3.9, 3.10]{kpr} and is
proved similarly.

\begin{thm} 
Let $\Lambda$ be a $k$-graph and suppose that the countable
group $G$ acts freely on $\Lambda$, then
\[
C^* ( \Lambda ) \rtimes_\beta G \cong 
C^* ( \Lambda / G ) \otimes {\mathcal K}\left( \ell^2 ( G ) \right) .
\]
Equivalently, if $c : \Lambda' \rightarrow G$ is a functor, then
\[
C^* ( G \times_c \Lambda' )\rtimes_\beta G\cong 
C^* ( \Lambda' ) \otimes {\mathcal K}\left ( \ell^2 ( G ) \right)
\]
where $\beta$, the action of $G$ on $C^* ( G \times_c \Lambda' )$, 
is induced by the natural action on $G \times_c \Lambda'$. If $G$ is
abelian this action is dual to $\alpha^c$ under the identification of
\ref{after}.
\end{thm}

\begin{proof}
The first statement follows from the second with $\Lambda' =
\Lambda /G$;  indeed,  by \ref{gtanalog} there is a functor $c :
\Lambda / G \rightarrow G$ such that 
$\Lambda \cong G \times_c (\Lambda / G )$ in an equivariant way.  The second
statement follows from applying \cite[Proposition 3.7]{kp} to the natural
$G$-action on
${\mathcal G}_{ G \times_c \Lambda'}  \cong 
{\mathcal G}_{\Lambda'} ( \tilde{c} )$. The final statement follows from
the identifications 
\[
C^* ( \Lambda ) \rtimes_{\alpha^c} \widehat{G} \cong  
C^* ( G \times_c \Lambda ) \cong  C^* ( {\mathcal G}_\Lambda ( \tilde{c} ) )
\]
and \cite[II.2.7]{r}.
\end{proof}
\section{$2$-graphs}

Given a $k$-graph $\Lambda$ one obtains for each $n \in {\bf N}^k$  a matrix 
$$M_\Lambda^n(u, v) = 
\# \{\lambda \in \Lambda^n : r(\lambda ) = u, s(\lambda ) = v \}.$$
By our standing assumption the entries are all finite and there are no
zero rows.   Note that for any $m, n \in  {\bf N}^k$ we have 
$M_\Lambda^{m+n} = M_\Lambda^mM_\Lambda^n$ (by the factorization property);
consequently, the matrices $M_\Lambda^m$ and $M_\Lambda^n$ commute for all
$m,n \in  {\bf N}^k$.   If $W$ is the $k$-graph associated to the
commuting  matrices  
$\{M_1 , \ldots , M_k \}$ satisfying conditions (H0)--(H3) of
\cite{rs2} which was 
considered in Example \ref{secondex}(iv), then one checks that 
$M_W^{e_i} = M_i^t$.
Further, if $\Lambda = E^*$ is a $1$-graph derived from the directed graph
$E$, then $M_\Lambda^1$ is the vertex matrix of $E$.  

Now suppose that  $A$ and $B$ are $1$-graphs with $A^0 = B^0 = V$ such the
associated vertex matrices commute.  Set 
$A^1*B^1 = \{ (\alpha, \beta) \in A^1 \times B^1 : s(\alpha) = r(\beta) \}$
and
$B^1*A^1 = \{ (\beta, \alpha) \in B^1 \times A^1 : s(\beta) = r(\alpha) \}$;
since the associated vertex matrices commute there is a bijection
$\theta: (\alpha, \beta) \mapsto (\beta', \alpha') $
from $A^1*B^1$ to $B^1*A^1$ such that $r(\alpha) = r(\beta')$ and 
$s(\beta) = s(\alpha')$.  We construct a $2$-graph $\Lambda$ from $A$,
$B$ and $\theta$.   This construction is very much in the spirit of 
\cite{rs2}; roughly speaking an element in $\Lambda$ of degree 
$(m,n) \in {\bf N}^2$ will consist of a rectangular grid of size
$(m,n)$ with edges of $A$ horizontally, edges of $B$ vertically and
nodes in $V$ arranged compatibly.
First identify $\Lambda^0  = V$.  For $(m,n) \in {\bf N}^2$ set 
$W(m,n) =  \{(i,j) \in  {\bf N}^2 :  (i,j) \le (m,n) \}$.   
An element in $\Lambda^{(m,n)}$ is given by
$v(i,j) \in V$ for $(i,j) \in W(m,n)$, $\alpha(i,j) \in A^1$ for 
$(i,j) \in W(m-1,n)$
and $\beta(i,j) \in B^1$ for $(i,j) \in W(m,n-1)$ (set $W(m,n) = \emptyset$
if $m$ or $n$ is negative) satisfying the following compatibility conditions
wherever they make sense:
\begin{itemize}
\item[i]   $r(\alpha(i,j) ) = v(i,j) $ and $r(\beta(i,j)) = v(i,j)$ 
\item[ii]   $s(\alpha(i,j) ) = v(i+1,j) $ and $s(\beta(i,j) ) = v(i,j+1) $
\item[iii] $\theta (\alpha (i,j),\beta(i+1,j)) = (\beta(i,j),\alpha(i,j+1))$;
\end{itemize}
for brevity and with a slight abuse of notation we regard this element
as a triple $(v, \alpha, \beta)$ (note that $\alpha$ disappears if 
$m = 0$ and $\beta$ disappears if $n=0$ and $v$ is determined by 
$\alpha$ and/or $\beta$ if $mn \neq 0$). Set 
$$\Lambda = \bigcup_{(m,n)} \Lambda^{(m,n)}$$ %% \in {\bf N}^2
and define $s(v, \alpha, \beta)  = v(m,n)$ and
$r(v,\alpha,\beta)  = v(0,0)$.

Note that if $\lambda \in A^m$ and $\mu \in B^n$ with $m,n > 0$ such that 
$s(\lambda) = r(\mu)$ there is a unique element  
$(v, \alpha, \beta) \in \Lambda^{(m,n)}$ such that 
$\lambda = \alpha(0,0) \alpha(1,0) \cdots\alpha(m-1,0)$ and
$\mu = \beta(m,0) \beta(m,1) \cdots \beta(m,n-1)$; denote this element 
$\lambda \mu$.
Further if $\lambda \in A^m$ and $\mu \in B^n$ with $m,n > 0$ such that 
$r(\lambda) = s(\mu)$ there is a unique element  
$(v, \alpha, \beta) $ in $\Lambda^{(m,n)}$ such that 
$\lambda = \alpha(0,n)\alpha(1,n)\cdots\alpha(m-1,n)$ and
$\mu = \beta(0,0) \beta(0,1) \cdots \beta(0,n-1)$; denote this element 
$\mu\lambda $.
Using these two facts it is not difficult to verify that given elements 
$(v, \alpha, \beta) \in \Lambda^{(m,n)}$  and
$(v', \alpha', \beta') \in \Lambda^{(m',n')}$ with $v(m, n) = v'(0, 0)$
there is a unique element
$(v'', \alpha'', \beta'') \in \Lambda^{(m+m', n+n')}$ such that
$v''(i,j) = v(i,j)$, $\alpha''(i,j) = \alpha(i,j)$, 
$\beta''(i,j) = \beta(i,j)$, $v''(m+i,n+j) = v'(i,j)$, 
$\alpha''(m+i,n+j) = \alpha'(i,j)$ and $\beta''(m+i,n+j) =
\beta'(i,j)$ wherever these formulas make sense.  Write 
$(v'', \alpha'', \beta'') = (v, \alpha, \beta) (v', \alpha', \beta')$.
This defines composition in $\Lambda$; note that associativity and the 
factorization property are built into the construction (as in \cite{rs2}).  
Finally, we write $\Lambda = A*_{\theta}B$.
It is straightforward to verify that up to isomorphism any $2$-graph
may be obtained from its constituent $1$-graphs in this way.

If $A = B$, then we may take $\theta = \iota$ the identity map.  In
that case one has $A*_{\iota}A \cong f^*(A)$ where 
$f: {\bf N}^2 \rightarrow \bf N$ is given by $f(m,n) = m+n$.   
Hence, by Corollary \ref{giutconsequences}(iii) we have
$C^* (A*_{\iota}A) \cong C^* (A)  \otimes C ( {\bf T} )$.

To further emphasise the dependence of the product $A *_{\theta} B$ on
the bijection $\theta : A^1*B^1 \rightarrow B^1*A^1$ consider the following
example: 

\begin{example}
Let $A = B $ be the $1$-graph derived from the directed graph which consists of one
vertex and two edges, say  $A^1 = \{ e , f \}$ (note $C^* (A) \cong {\mathcal O}_2$).  Then
$A^1*A^1 = \{ (e, e), (e, f), (f, e),  (f, f) \}$, and we define the bijection $\theta$ to
be the flip.  It is easy to show that 
$A *_{\theta} A \cong A \times A$; hence,  
$$ 
C^* ( A *_{\theta} A ) \cong 
{\mathcal O}_2 \otimes {\mathcal O}_2 \cong {\mathcal O}_2 
$$
where the first isomorphism follows from Corollary
\ref{giutconsequences}(iv) and the second from the Kirchberg-Phillips 
classification theorem (see \cite{k,p}).  But  
$$
C^*(A *_{\iota} A ) \cong {\mathcal O}_2 \otimes C ( {\bf T} );
$$
hence,  $A *_{\theta} A \not\cong A *_{\iota} A$.
\end{example}


\begin{thebibliography}{KPRR}
\bibitem[A-D]{a-d} C.~Anantharaman--Delaroche.
\newblock Purely infinite $C^*$-algebras arising from dynamical
systems.
\newblock {\em Bull. Soc. Math. France}, {\bf 125}:  199--225,
(1997).

\bibitem[A-DR]{a-dr}C.~Anantharaman--Delaroche and J.~Renault.
\newblock Amenable groupoids. {\em To appear.}
%\newblock{\em Memoirs of Amer. Math. Soc.}

\bibitem[ALNR]{alnr}
S. Adji, M. Laca, M. Nilsen and I. Raeburn.
\newblock Crossed products by semigroups of endomorphisms and
the Toeplitz algebras of ordered groups.
\newblock {\em Proc. Amer Math. Soc.}, {\bf 122}: 1133--1141,
(1994).

\bibitem[BPRS]{bprs}
T.~Bates, D.~Pask, I.~Raeburn, W.~Szymanski.
\newblock The $C^*$--algebras of row--finite graphs.
\newblock {\em Submitted}.

%\bibitem[BKR]{bkr} S.~Boyd, N.~Keswani,~I.~Raeburn.
%\newblock Faithful representations of crossed products by endomorphisms.
%\newblock {\em Proc. Amer. Math. Soc.}, {\bf 118}: 427--436, (1993).

\bibitem[B]{b}
O.~Bratteli.
\newblock Inductive limits of finite dimensional $C^*$--algebras.
\newblock {\em Trans. Amer. Math. Soc.}, {\bf 171}: 195--234,
(1972).

\bibitem[CK]{ck}
J.~Cuntz and W.~Krieger.
\newblock A class of ${C}^*$-algebras and topological
{M}arkov  chains.
\newblock {\em Invent. Math.}, {\bf 56}: 251--268, (1980).

\bibitem[D]{d} V.~Deaconu.
\newblock Groupoids associated with endomorphisms
\newblock{\em Trans. Amer. Math. Soc.}, {\bf 347}: 1779--1786,
(1995).

\bibitem[GT]{gt} J.L.~Gross and T.W.~Tucker.
\newblock {\em Topological graph theory.}
\newblock Wiley Interscience Series in Discrete
Mathematics and Optimization, First edition
(1987)

\bibitem[H]{h} P.J. Higgins.
\newblock {\em Notes on categories and groupoids.}
\newblock van Nostrand Rienhold (1971).

\bibitem[aHR]{ahr}
A. ~an Huef and I.~Raeburn.
\newblock The ideal structure of {C}untz-{K}rieger algebras.
\newblock {\em Ergod. Th. and Dyn. Sys.}, {\bf 17}: 611--624,
(1997).

\bibitem[KQR]{kqr}
S.~Kaliszewski, J.~Quigg and I.~Raeburn.
\newblock  Skew products and crossed products.
by coactions.
\newblock {\em Preprint}.

\bibitem[K]{k}  E.~Kirchberg, 
\newblock The classification of purely infinite $C^*$-algebras
using Kasparov's theory.
\newblock {\em Preprint}.

\bibitem[KPRR]{kprr} 
A.~Kumjian, D.~Pask,  I.~Raeburn, and J.~Renault.
\newblock Graphs, groupoids and Cuntz--Krieger algebras.
\newblock{\em J. Funct. Anal.}, {\bf 144}: 505--541, (1997).

\bibitem[KPR]{kpr} A.~Kumjian, D.~Pask, I.~Raeburn.
\newblock Cuntz--Krieger algebras of directed graphs,
\newblock {\em Pacific. J. Math.}, {\bf 184}: 161--174, (1998).

\bibitem[KP]{kp} A.~Kumjian and D.~Pask.
\newblock $C^*$--algebras of directed graphs and group actions,
\newblock {\em Ergod. Th. \& Dyn. Sys.}, to appear.

\bibitem[LS]{ls} M.~Laca and J.~Spielberg.
\newblock Purely infinite $C^*$--algebras from boundary actions of
discrete groups,
\newblock {\em J. Reine Angew. Math.}, {\bf 480}: 125--139, (1996).

\bibitem[MacL]{macl} S.~MacLane.
\newblock Categories for the working Mathematician,
\newblock Graduate Texts in Mathematics {\bf 5}, Springer--Verlag, 1971.

\bibitem[Ma]{ma} T.~Masuda.
\newblock Groupoid dynamical systems and crossed product II --
the case of ${C}^*$--systems.
\newblock {\em Publ. RIMS Kyoto Univ.}, {\bf 20}: 959--970,
(1984).

\bibitem[Mu]{mu} P.~Muhly.
\newblock A finite dimensional introduction to Operator algebra,
\newblock In {\em Operator algebras and applications (Samos, 1996)},
313--354, NATO Adv. Sci. Inst. Ser. C Math. Phys. Sci., 495, Kluwer
Acad.  Publ., Dordrecht, 1997.

\bibitem[P]{p} N.C. Phillips,
\newblock  A classification theorem for nuclear purely infinite simple
$C^*$-algebras.
\newblock {\em Preprint}.

\bibitem[R]{r} J.~Renault.
\newblock { A groupoid approach to ${C}^*$-algebras}. 
{\em Lecture Notes in Mathematics}, vol. ~{\bf 793}.
\newblock Springer-Verlag, 1980.

\bibitem[RS1]{rs1} G.~Robertson and T.~Steger.
\newblock  $C^*$--algebras arising from group actions on the
boundary of a triangle building,
\newblock {\em Proc. London Math. Soc.}, {\bf 72}: 613--637, (1996).

\bibitem[RS2]{rs2} G.~Robertson and T.~Steger.
\newblock Affine buildings, tiling systems and higher rank
Cuntz--Krieger algebras,
\newblock {\em J. Reine Angew. Math.},  {\bf 513}: 115--144, (1999).

\bibitem[RS3]{rs3} G.~Robertson and T.~Steger.
\newblock $K$--theory for rank two Cuntz--Krieger algebras.
\newblock {\em Preprint}.

\bibitem[RSc]{rs-uct} J.~Rosenberg and C.~Schochet.
\newblock The {K}\"{u}nneth theorem and the universal
coefficient theorem for {K}asparov's generalized $K$--functor.
\newblock {\em Duke Math. J.}, {\bf 55}: 431--474, (1987).

\end{thebibliography}
\end{document}